\documentclass[11pt,twoside]{article}

\usepackage{mathrsfs,amssymb,amsfonts,amsthm,amsmath,amsbsy}
\usepackage{graphicx}
\graphicspath{{figures/}}
\usepackage{dsfont}
\usepackage[dvipsnames]{xcolor}
\usepackage[authoryear]{natbib}
\usepackage[bookmarks,bookmarksnumbered,colorlinks=true,pdfstartview=FitV, linkcolor=Red,citecolor=blue!90!black,urlcolor=blue!90!black]{hyperref}
\hypersetup{
  pdftitle={A tractable transition function for the Wright--Fisher diffusion with selection via gamma duality},
  pdfauthor={Jaromir Sant, Paul A. Jenkins, Matteo Ruggiero}
}
\usepackage{autonum}
\usepackage{cancel}
\usepackage{leftindex}
\usepackage[ruled, lined, linesnumbered]{algorithm2e}
\usepackage{fancyhdr}
\renewcommand{\baselinestretch}{1.1}
\allowdisplaybreaks

\bibpunct{\textcolor{blue!90!black}{(}}{\textcolor{blue!90!black}{)}}{\textcolor{blue!90!black}{;}}{a}{\textcolor{blue!90!black}{,}}{\textcolor{blue!90!black}{;}}

\def\XX{\boldsymbol{X}}
\def\xx{\boldsymbol{x}}
\def\YY{\boldsymbol{Y}}
\def\yy{\boldsymbol{y}}
\def\tt{\boldsymbol{t}}
\def\mm{\boldsymbol{m}}
\def\lv{\boldsymbol{\ell}}
\def\LL{\boldsymbol{L}}

\def\nn{\boldsymbol{n}}
\def\NN{\boldsymbol{N}}
\def\kk{\boldsymbol{k}}
\def\bb{\boldsymbol{b}}
\def\oo{\boldsymbol{0}}
\def\aa{\boldsymbol{\alpha}}
\def\ss{\boldsymbol{\sigma}}
\def\MN{\mathrm{MN}}
\def\PP{\mathbb{P}}
\def\ZZ{\mathbb{Z}^{K}_{+}}
\def\D{\mathcal{D}}
\def\A{\mathcal{A}}

\numberwithin{equation}{section}
\newtheorem{theorem}{Theorem}[section]
\newtheorem{proposition}[theorem]{Proposition} 

\newtheorem{lemma}[theorem]{Lemma} 
 
\newtheorem{definition1}[theorem]{Definition} \newenvironment{definition}{\begin{definition1}\rm}{\hfill$\square$\end{definition1}}
\newtheorem{remark1}[theorem]{Remark} \newenvironment{remark}{\begin{remark1}\rm}{\hfill$\square$\end{remark1}}

\long\def\symbolfootnote[#1]#2{\begingroup\def\thefootnote{\hspace*{-1mm}\fnsymbol{footnote}}\footnote[#1]{#2}\endgroup}

\fancypagestyle{plain}{
\fancyhead{}
}

\begin{document}

\title{A tractable transition function for the Wright--Fisher diffusion with selection via gamma duality}
\author{Jaromir Sant$^{1}$, Paul A. Jenkins$^{2}$ and Matteo Ruggiero$^{3}$\\[2mm]
\small $^{1}$Universidad Carlos III de Madrid\\
\small $^{2}$University of Warwick\\
\small $^{3}$Stern School of Business, NYU Abu Dhabi}
\date{}
\maketitle

\begin{abstract}
The transition function of the Wright--Fisher diffusion with selection is a fundamental object for non-neutral evolutionary modelling and likelihood-based inference, but, unlike in the neutral case, existing representations are impractical or inefficient. Spectral approaches rely on truncation levels that grow with selection strength and the number of alleles, thus leading to rapidly increasing computational costs, while existing duality-based expansions lead to multi-type dual processes whose transition probabilities involve intractable normalising constants.

In this paper we obtain a new transition-function representation for the multi-allelic Wright--Fisher diffusion with parent-independent mutation and genic selection that makes transition-law evaluation and simulation practically feasible in regimes where existing methods are costly or unavailable. The key tool is a new \emph{gamma duality}, obtained by augmenting classical polynomial duality by an exponential tilting and yields a multivariate birth-death dual with fully explicit infinitesimal rates and an entrance-boundary construction. This is used to express the Wright--Fisher transition law as an explicit mixture of \emph{standard} Dirichlet components, with mixing weights obtained from the dual started from its entrance boundary. We characterise the dual stationary distribution and its small-time behaviour, and develop practical and efficient sampling procedures for both transition and stationary laws, including a fast approximate stationary sampler with explicit error control. Numerical experiments illustrate substantial runtime gains over specialised methods in regimes where they apply, while our approach remains applicable for arbitrary numbers of alleles and selection coefficients.
\end{abstract}

\renewcommand{\baselinestretch}{1.1}\selectfont

\section{Introduction}\label{Introduction}

The Wright--Fisher (WF) diffusion arises as a scaling limit of finite-population Markov chains, including the WF and Moran models and, more generally, Cannings-type models. It provides a canonical model for the evolution of allele frequencies in a population under random genetic drift, mutation, and selection \citep[Ch.\ 10]{EthierKurtz}, and also appears in a range of other applied fields, including econometrics and mathematical finance \citep{Gourieroux,Larsen,Filipovic,Arrouy}, game theory \citep{chen2018wright}, topic modelling \citep{Perrone}, Bayesian nonparametric statistics \citep{mena2016dynamic}, opinion formation \citep{Chen,Furioli}, and models of neuronal activity \citep{DOnofrio2}. See also \cite{pal2013wright,costantini2017wright,garcia2021exact,gonzalez2018duality,cordero2025lambda} for recent work extending the classical WF diffusion.

The primary domain of application is in population genetics, where the relative importance of drift, mutation, and selection in driving evolution has long been a central question \citep{ewe:2004:I}. This focus naturally leads to inferring evolutionary forces from discretely observed allele-frequency time series, for which the WF diffusion has become a standard workhorse \citep{Bollback,Steinrucken2014novel,Schraiber2016bayesian,He,tataru2017statistical,buzbas2009maximum}. Parameter estimation in such settings requires the repeated evaluation of, or sampling from, transition laws across multiple time points and parameter values, both within likelihood-based inference (including diffusion-driven hidden Markov models), as well as in simulation-based approaches such as exact data augmentation schemes and particle filters. Therefore, the absence of a tractable transition function representation (and consequently of an efficient sampler with controlled error) constitutes a key computational bottleneck in this regard.

For the neutral WF diffusion (i.e., in the absence of selection), the identification of the transition
function dates back to Kimura's classical \emph{spectral} expansions, summarized in \cite{kim:1964}, which express the transition density in terms of an orthogonal polynomial family of eigenfunctions, with \cite{gri:1979:AAP11:310} extending these results to incorporate recurrent mutation. Later, a \emph{probabilistic} series expansion was derived, expressing the transition density as an infinite mixture of Dirichlet distributions \citep{Griffiths80,gri:li:1983,tavare1984line}. See also \cite{feng2011functional,griffiths2025dual} for analogous representations in closely related infinite-dimensional settings. The spectral and probabilistic representations are equivalent \citep{GS2010}, but the Dirichlet-mixture form is particularly appealing computationally, as it enables exact simulation routines for neutral WF paths and bridges \citep{gri:li:1983,JenkinsSpano,SantEWF,sant2025ewf}, thereby enabling likelihood-based inference for discretely observed diffusions and diffusion-driven hidden Markov models.


By contrast, obtaining a workable transition-function representation for the non-neutral WF diffusion is substantially more challenging. Classical spectral approaches lead to eigenfunctions which are typically available only through special function representations instead of multivariate polynomials \citep{kim:1955:CSHSQB}. Probabilistic arguments still yield an exact mixture expansion in terms of exponentially tilted Dirichlet components \citep{BarbourEthierGriffiths}, but the associated multi-type birth-death dual has state-dependent rates involving intractable normalising constants, a difficulty that propagates to likelihood evaluation and reduces its practical applicability \citep{donnelly2001likelihoods,buzbas2009maximum}. Alternatively, one can suitably truncate the spectral representation obtained by considering an exponentially tilted version of the neutral basis eigenfunctions \citep{song2012simple,steinrucken2013explicit}, which enables likelihood methods \citep{Steinrucken2014novel} to be used (including the diploid selection case). 
However, the computational cost grows rapidly with the number of alleles and especially with the selection strength, as the truncation level required for accuracy increases accordingly. 
Exact transition-law simulation via retrospective rejection sampling provides a viable alternative \citep{schraiber2013analysis,JenkinsSpano}, but acceptance rates deteriorate as either the magnitude of selection or the time horizon increases, with current implementations being available only in the one-dimensional setting \citep{SantEWF}.

In this paper we obtain a new transition-function representation for the Wright--Fisher diffusion with parent-independent mutation and genic selection that makes transition-law evaluation and simulation practically feasible in regimes where existing methods are costly or unavailable. The key tool enabling this representation is a new \emph{gamma duality} for the WF diffusion with selection. Duality \citep{jansen2014notion} is pervasive in WF-type models and links allele-frequency diffusions to genealogical or particle systems, yielding moment identities and sampling representations \citep{EthierKurtz,M99,EG09,Eea10,Cea15,gonzalez2018duality} and lending itself to inference for certain hidden Markov models \citep{papas_ruggiero,ascolani2023smoothing,king2024approximate,dallapria2026}. In particular, \citet{BarbourEthierGriffiths} exploit a polynomial moment duality to derive a transition-function expansion in the non-neutral case, but the resulting dual dynamics involve ratios of normalising constants that are not tractable to evaluate or simulate from. In our formulation, gamma duality augments classical polynomial duality by an exponential tilting and yields a multivariate birth-death dual with fully explicit infinitesimal rates and an entrance-boundary construction. This dual turns the intractable mixture of \citet{BarbourEthierGriffiths} into a mixture of \emph{standard} Dirichlet components, with mixing weights obtained from the dual started from its entrance boundary. Since the dual dynamics are straightforward to simulate, this removes the normalising-constant bottleneck and allows us to build practical and efficient samplers for both transition and stationary laws.

Relative to the existing approaches, our method applies for arbitrary numbers of alleles and across a wide range of selection coefficients, including regimes where retrospective rejection or spectral truncations are unavailable or impractical. Our dual-based samplers improve on existing strategies in scope, computation time, or both, depending on the regime.

The remainder of the paper is organised as follows. Section~\ref{sec:model} introduces the model, notation, and background on existing transition function representations and their limitations. Section~\ref{GammaDuality} develops gamma duality and characterises the corresponding dual. Section~\ref{sec:transfunc} presents the resulting Dirichlet-mixture representation of the non-neutral transition function, along with sampling algorithms for the transition and stationary laws, and accompanying numerical illustrations. A brief discussion is given in Section~\ref{Discussion}, while additional numerical experiments briefly mentioned in the illustration are provided in the Appendix.


\section{The model and current limitations}\label{sec:model}

\subsection{The WF diffusion with mutation and genic selection}\label{subsec:model}

Denote by $\XX := (\XX_t)_{t\geq0}$ the $K$-allele WF diffusion with genic selection and parent-independent mutation, taking values in the $(K-1)$-dimensional simplex
\[
\Delta_{K} := \Big\{ \xx\in[0,1]^{K} : \sum_{i=1}^{K} x_i = 1 \Big\}.
\]
Its infinitesimal generator is 
\begin{align}\label{eq:WFGeneratorHaploidPIM}
	\A_K
	=
	\frac{1}{2}\sum_{i,j=1}^{K}x_i\left(\delta_{ij} - x_j\right)\frac{\partial^2}{\partial x_i \partial x_j}
	\!+\! \frac{1}{2}\sum_{i=1}^{K}\left[\alpha_i - \left|\aa\right|x_i + x_i\left(\tilde\sigma_i - \sum_{j=1}^{K}\tilde\sigma_jx_j\right)\right]\frac{\partial}{\partial x_i},
\end{align}
where $\tilde\ss = (\tilde\sigma_{i})_{i=1}^{K} \in \mathbb{R}^{K}$ describes the selection pattern, $\aa = (\alpha_{i})_{i=1}^{K} \in [0,\infty)^K$ are the mutation parameters, and $\left|\aa\right| := \sum_{i=1}^{K}\alpha_i$. 
We take the domain of $\A_K$ to be $\D(\A_{K}) := \{ f|_{\Delta_K} : f \in C^{2}(\mathbb{R}^{K}) \}$, where $f|_{\Delta_K}$ denotes restriction to $\Delta_K$. 

The model is invariant to adding a constant to all selection coefficients, as replacing $\tilde\ss$ by $\tilde\ss+c\boldsymbol{1}$, with $\boldsymbol{1} = (1)_{i=1}^{K}$, leaves the drift in \eqref{eq:WFGeneratorHaploidPIM} unchanged.
Accordingly, without loss of generality we work henceforth with the shifted coefficients $\ss = (\sigma_{i})_{i=1}^{K} \in [0,\infty)^{K}$, where
\begin{align}\label{eq:RescalingSigma_sec2}
	\sigma_i := \tilde\sigma_i - \min_{k \in \{1,\dots,K\}}\tilde\sigma_k \ge 0, \qquad i=1,\dots,K.
\end{align}
When $\alpha_i > 0$ for all $i$, the stationary distribution $\Pi^{\ss}_{\aa}$ of the diffusion is unique and has a density with respect to Lebesgue measure on $\Delta_K$ given by
\begin{align}\label{eq:StationaryWF}
	\pi^{\ss}_{\aa}(\xx)
	&=
	c_{\aa, \ss}^{-1}e^{\left<\ss, \xx\right>}\xx^{\aa - \boldsymbol{1}},
	\qquad
	c_{\aa, \ss} := \int_{\Delta_K}e^{\left<\ss, \xx\right>}\xx^{\aa - \boldsymbol{1}}d\xx,
\end{align}
where $\xx^{\aa - \boldsymbol{1}} := \prod_{i=1}^K x_i^{\alpha_i - 1}$ and $\left< \ss, \xx\right> := \sum_{i=1}^{K}\sigma_ix_i$.
See \citet[Ch.\ 10]{EthierKurtz} for further details.

Setting $\ss = \oo$ in \eqref{eq:WFGeneratorHaploidPIM} yields the generator of the neutral WF diffusion with parent-independent mutation. Its stationary density, obtained from \eqref{eq:StationaryWF}, is the standard Dirichlet density $\pi^{\oo}_{\aa}(\xx)$ with normalising constant
\begin{equation}\nonumber
	c_{\aa, \oo}
	=
	\frac{\prod_{j=1}^{K}\Gamma(\alpha_j)}{\Gamma(\left|\aa\right|)}.
\end{equation} 
The transition function of the neutral WF diffusion admits the Dirichlet-mixture representation \citep{Griffiths80,tavare1984line}
\begin{align}\label{eq:NeutralWFTransitionDensity}
	P_{\aa}^{\oo}(\xx, d\yy; t)
	&=
	\sum_{m=0}^{\infty} d_{m}(t)\sum_{\lv\in\ZZ : \left|\lv\right| = m}\MN(\lv;\, m, \xx)\,\Pi^{\oo}_{\aa+\lv}(d\yy),
\end{align}
where $\mathbb{Z}_{+}:=\{0,1,\ldots\}$, $\MN(\lv;\, m, \xx) = \binom{m}{\lv}\xx^{\lv}$ is the multinomial pmf with parameters $(m,\xx)$ and $\Pi^{\oo}_{\aa+\lv}$ is the Dirichlet distribution with parameter $\aa+\lv$ and density $\pi^{\oo}_{\aa+\lv}(\cdot)$ obtained from \eqref{eq:StationaryWF} when $\ss = \oo$; cf.\ \citet[Cor.~1.9]{ethier1993transition}.
The mixing weights $d_m(t)$ form a probability mass function on $\mathbb{Z}_+$ and have a standard genealogical interpretation via Kingman's coalescent \citep{Kingman82a}: they coincide with the time-$t$ distribution of the one-dimensional pure-death block-counting process where lineages are lost by both coalescence and mutation.
The representation \eqref{eq:NeutralWFTransitionDensity} underpins exact simulation of the neutral WF diffusion and bridge \citep{JenkinsSpano}, avoiding time-discretization schemes that perform poorly on bounded domains, or truncations whose error is difficult to quantify.

For later comparison, note that \eqref{eq:NeutralWFTransitionDensity} can equivalently be written as a single mixture over $\lv\in\ZZ$, i.e.,
\begin{align}\label{eq:NeutralWFTransitionDensity_single}
	P_{\aa}^{\oo}(\xx,d\yy;t)
	&=
	\sum_{\lv\in\ZZ}
	b^{\oo}_{\lv}(t,\xx)\,\Pi^{\oo}_{\aa+\lv}(d\yy),
	\qquad
	b^{\oo}_{\lv}(t,\xx)
	:=
	d_{|\lv|}(t)\,\MN(\lv;\, |\lv|,\xx),
\end{align}
where the neutral weights factor into the law of the block count $|\lv|$ and a multinomial type configuration conditional on $|\lv|$. Further note that the representations \eqref{eq:NeutralWFTransitionDensity} and \eqref{eq:NeutralWFTransitionDensity_single} remain valid when some of the $\alpha_i$ are zero provided we define the Dirichlet($\aa$) distribution via
\begin{equation} \label{eq:Dirichlet-generalization}
	\Pi^{\oo}_{\aa}(\cdot) := \mathbb{P}\left[(Y_1/Y,\dots,Y_K/Y) \in \cdot\right],
\end{equation}
where $Y_1,\dots,Y_K$ are independent random variables with $Y_i \sim\Gamma(\alpha_i,1)$, $Y = \sum_{i=1}^K Y_i$, and we recall that $\mathbb{P}(Y_i = 0) = 1$ if $Y_i \sim \Gamma(0,1)$ \citep[eq.~(1.26)]{ethier1993transition}. More generally, $\Pi^{\ss}_{\aa}$ is defined via
\begin{equation}
	\label{eq:StationaryWF-RNderivative}
	\frac{d\Pi^{\ss}_{\aa}}{d\Pi^{\oo}_{\aa}}(\yy) = \frac{c_{\aa,\oo}}{c_{\aa, \ss}}e^{\left<\ss, \yy\right>}.
\end{equation}
When $\aa = \oo$, $\Pi^{\oo}_{\oo}$ is undefined, but in that case $d_0(t) = 0$ and the component $\Pi^{\oo}_{\oo}$ does not appear in \eqref{eq:NeutralWFTransitionDensity}--\eqref{eq:NeutralWFTransitionDensity_single}. Our results will also allow for the possibility that some or all $\alpha_i = 0$, and so we take \eqref{eq:Dirichlet-generalization}--\eqref{eq:StationaryWF-RNderivative} as the definition for $\Pi^{\ss}_{\aa}$ hereonwards.


\subsection{Existing transition function representations and their limitations}\label{subsec:limitations}

For $\ss\neq\oo$, transition-function representations analogous in spirit to
\eqref{eq:NeutralWFTransitionDensity_single} remain available, but the ingredients that make the neutral case computationally convenient typically cease to be explicit or easy to evaluate in regimes of practical interest. Two broad computational routes are commonly used. One utilises spectral theory: eigenfunctions can be expanded in an exponentially tilted neutral (Jacobi) basis and truncated for numerical evaluation \citep{song2012simple,steinrucken2013explicit}. This enables the application of likelihood methods under diploid selection \citep{Steinrucken2014novel}, at the price of dealing with truncation error control and computational costs that increase with selection strength and dimension. The alternative, probabilistic approach is based on duality, which is the focus here.

Duality is often the mechanism that yields usable transition expansions. In the neutral setting, polynomial (moment) duality links the WF diffusion to Kingman's coalescent and underpins classical mixture representations of the transition function \citep{ethier1993transition,EG09,GS2010}. 
For general diploid selection, \cite{BarbourEthierGriffiths} derive a duality-based transition expansion; in the case of genic selection, this reduces to a multi-type birth-death dual on $\ZZ$ with rates 
\begin{equation}\label{eq:BEG_rates_sec2}
q(\nn, \nn - \boldsymbol{e}_i)=
\frac{n_i\left(|\aa| + |\nn| - 1\right)}{2}\,
\frac{c_{\aa + \nn - \boldsymbol{e}_i, \ss}}{c_{\aa + \nn, \ss}},\quad 
q(\nn, \nn + \boldsymbol{e}_i)=
|\nn|\,\hat\sigma_i\,\frac{n_i + \alpha_i}{|\nn| + |\aa|}\,
\frac{c_{\aa + \nn + \boldsymbol{e}_i, \ss}}{c_{\aa + \nn, \ss}},
\end{equation}
where $\hat\sigma_{i}=\max_{j}\tilde{\sigma}_{j}-\tilde{\sigma}_{i}$, $\boldsymbol{e}_i$ is the $i^{\textnormal{th}}$ canonical basis vector, and the ratios involve the normalising constants from \eqref{eq:StationaryWF}.
Combining this duality with sample inversion yields the transition expansion
\begin{align}\label{eq:BEG_transition_expansion}
P_{\aa}^{\ss}(\xx,d\yy;t)
&=
\sum_{\lv\in\ZZ} b^{\ss}_{\lv}(t;\xx)\,\Pi^{\ss}_{\aa+\lv}(d\yy),
\end{align}
where $\Pi^{\ss}_{\aa+\lv}$ is defined in \eqref{eq:StationaryWF-RNderivative} and the mixing weights
$b^{\ss}_{\lv}(t;\xx)$ are characterised by the birth-death process with rates
\eqref{eq:BEG_rates_sec2} started from its entrance boundary. Unlike the neutral case, the weights do not admit a simple closed form nor a factorisation into a one-dimensional law for $|\lv|$ and a conditional multinomial split, because selection couples types through the state-dependent normalising constants.

While \eqref{eq:BEG_transition_expansion} is analytically informative, it has limited practical use: simulating the dual and evaluating or approximating $b^{\ss}_{\lv}(t;\xx)$ requires repeated computation of ratios of $c_{\aa+\nn,\ss}$. Numerical methods for evaluating $c_{\aa,\ss}$ (and related exponentially weighted Dirichlet integrals) are available in structured settings \citep{genz2003computation,joyce2012efficient,fearnhead2004exact}, but they do not remove the core bottleneck here: repeated evaluation of such ratios quickly becomes expensive as the dimension $K$ and the strength of selection grow.

The next section introduces a modified duality that preserves the conceptual advantages of \cite{BarbourEthierGriffiths} while yielding a fully tractable dual process, leading to a Dirichlet-mixture representation with standard Dirichlet components and computationally accessible mixing weights.


\section{Gamma duality}\label{GammaDuality}

\subsection{The dual process}

We establish a new duality for the multi-allelic WF diffusion on $\Delta_{K}$ with parent-independent mutation and genic selection, obtained by an exponential tilting of the classical polynomial duality functions $\xx^{\nn}/m(\nn)$, where $m(\cdot)$ is some suitably chosen function.

Consider the infinitesimal generator $\A_K$ defined in \eqref{eq:WFGeneratorHaploidPIM} with $\aa \in (0,\infty)^K$, and define the function $h : \Delta_K \times \ZZ \to [0,\infty)$ by
\begin{align}\label{eq:GammaDualityFunction}
	h(\xx, \nn)
	=
	\frac{c_{\aa, \ss}}{c_{\aa+\nn, \oo}}
	\,e^{-\langle\ss, \xx\rangle}\,\xx^{\nn}
	=: \frac{c_{\aa, \ss}}{c_{\aa+\nn, \oo}}\tilde{h}(\xx, \nn),
\end{align}
with constants as in \eqref{eq:StationaryWF}.

Consider also the following process.
\begin{definition}\label{def:Dual}
	Let $(\NN_t)_{t\ge0}$ be a birth-and-death process on $\ZZ$ with transition rates
	\begin{align}\label{eq:TransitionRatesDual}
		q(\nn, \nn - \boldsymbol{e}_i) = 
		\frac{n_i\big(\,|\aa| + |\nn| - 1\,\big)}{2},
		\quad 
		q(\nn, \nn + \boldsymbol{e}_i) =
		\frac{\sigma_i\big(\alpha_i + n_i\big)}{2}, \quad 
		i=1,\dots,K,
	\end{align}
	and $q(\nn,\nn')=0$ for all other $\nn'\neq\nn$, where $\aa \in [0,\infty)^K$, and $\sigma_i$ is as in \eqref{eq:RescalingSigma_sec2}.
\end{definition}

The following theorem establishes the duality between the non-neutral WF diffusion and the birth-and-death process in Definition \ref{def:Dual}.

\begin{proposition}\label{thm:Duality}
	Let $(\NN_t)_{t\ge0}$ be as in Definition \ref{def:Dual} but with $\aa\in(0,\infty)^{K}$. Then, for all $\xx\in\Delta_K$, $\nn\in\ZZ$ and $t\ge0$,
	\begin{align}\label{eq:GammaDuality}
		\mathbb{E}\!\left[h\!\left(\XX_t, \nn\right) \,\middle|\, \XX_0 = \xx \right]
		=
		\mathbb{E}\!\left[h\!\left(\xx, \NN_t\right) \,\middle|\, \NN_0 = \nn\right],
	\end{align}
	where $(\XX_t)_{t\ge0}$ is the WF diffusion with generator \eqref{eq:WFGeneratorHaploidPIM}, and $h$ is defined in \eqref{eq:GammaDualityFunction}.
\end{proposition}

\begin{proof}
	It is clear that $h(\cdot, \nn) \in \D(\A_{K})$, and that
\begin{align}\nonumber
			\A_K \tilde{h}(\xx, \nn) &= \frac{1}{2}\sum_{i,j=1}^{K}x_i\left(\delta_{ij} - x_j\right)\Big[n_i\left[ \left(n_j - \delta_{ij}\right)\tilde{h}\left(\xx, \nn - \boldsymbol{e}_i - \boldsymbol{e}_j\right) - \tilde{\sigma}_j \tilde{h}\left(\xx, \nn - \boldsymbol{e}_i\right)\right]  \\
			&{}\qquad\qquad\qquad\qquad\qquad\qquad\quad{}- \tilde{\sigma}_i\left[n_j \tilde{h}\left(\xx, \nn - \boldsymbol{e}_j\right) - \tilde{\sigma}_j\tilde{h}\left(\xx, \nn\right)\right]\Big] \\
			&{}\quad + \frac{1}{2}\sum_{i=1}^{K}\left[\alpha_i - \left|\aa\right|x_i + x_i\left(\tilde{\sigma}_i - \sum_{j=1}^{K}\tilde{\sigma}_jx_j\right)\right]\left[n_i \tilde{h}\left(\xx, \nn-e_i\right) - \tilde{\sigma}_i\tilde{h}\left(\xx, \nn\right)\right]  \\
			&{}= \frac{1}{2}\Bigg\{ \sum_{i=1}^{K}n_i\left(\alpha_i + n_i - 1\right)\tilde{h}\left(\xx, \nn- \boldsymbol{e}_i\right) - \left|\nn\right|\left(\left|\aa\right| + \left|\nn\right| - 1\right)\tilde{h}\left(\xx, \nn\right)  \\
			&\qquad\qquad\qquad{} + \sum_{i=1}^{K}\tilde{\sigma}_i\left(\left|\aa\right| + \left|\nn\right|\right)\tilde{h}\left(\xx, \nn + \boldsymbol{e}_i\right) - \sum_{i=1}^{K}\tilde{\sigma}_i\left(\alpha_i + n_i\right)\tilde{h}\left(\xx, \nn\right)\Bigg\}  \\
			&{}= \frac{1}{2}\Bigg\{ \sum_{i=1}^{K}n_i\left(\alpha_i + n_i - 1\right)\tilde{h}\left(\xx, \nn- \boldsymbol{e}_i\right) - \left|\nn\right|\left(\left|\aa\right| + \left|\nn\right| - 1\right)\tilde{h}\left(\xx, \nn\right)  \\
			&\qquad\qquad\qquad{} + \sum_{i=1}^{K}\sigma_i\left(\left|\aa\right| + \left|\nn\right|\right)\tilde{h}\left(\xx, \nn + \boldsymbol{e}_i\right) - \sum_{i=1}^{K}\sigma_i\left(\alpha_i + n_i\right)\tilde{h}\left(\xx, \nn\right)\Bigg\}
\end{align}
	where we used \eqref{eq:RescalingSigma_sec2} in the third equality. Now since
	\begin{equation}\label{eq:TransRatesSimplif}
		\begin{aligned}
			\frac{\left(\alpha_i + n_i - 1\right)}{c_{\aa+\nn, \oo}} &= \frac{\Gamma\left(\left|\aa\right| + \left|\nn\right|\right)\left(\alpha_i + n_i - 1\right)}{\prod_{j=1}^{K}\Gamma\left(\alpha_j + n_j\right)} = \frac{\left|\aa\right| + \left|\nn\right| - 1}{c_{\aa + \nn - \boldsymbol{e}_i, \oo}},  \\
			\frac{\left|\aa\right| + \left|\nn\right|}{c_{\aa + \nn, \oo}} &= \frac{\Gamma\left(\left|\aa\right| + \left|\nn\right|\right)\left(\left|\aa\right| + \left|\nn\right|\right)}{\prod_{j=1}^{K}\Gamma\left(\alpha_j + n_j\right)} = \frac{\alpha_i + n_i}{c_{\aa + \nn + \boldsymbol{e}_i, \oo}},
		\end{aligned}
	\end{equation} 
	it follows that
	\begin{equation}\label{eq:DualGenerator}
		\begin{aligned}
			\A_K h\left(\xx, \nn\right) &= \frac{\left|\aa\right| + \left|\nn\right| - 1}{2}\sum_{i=1}^{K}n_i\left[h(\xx, \nn-\boldsymbol{e}_i) - h(\xx, \nn)\right]  \\
			&\quad{}+ \frac{1}{2}\sum_{i=1}^{K}\sigma_i\left(\alpha_i+ n_i\right)\left[h(\xx, \nn + \boldsymbol{e}_i) - h(\xx, \nn) \right].
		\end{aligned}
	\end{equation} 
	Lemma \ref{lemma:hInDomainGenerator} (see below) gives that $h(\xx,\cdot) \in\D(\A^{\#}_K)$, and thus \eqref{eq:DualGenerator} defines the generator of a birth-and-death process with transition rates \eqref{eq:TransitionRatesDual}, from which \citep[see Corollary 4.4.15]{EthierKurtz} the conclusion of the theorem follows. 
\end{proof}
\begin{lemma}\label{lemma:hInDomainGenerator}
	The function $h(\xx, \cdot)$ as defined in \eqref{eq:GammaDualityFunction} is in the domain of the generator $\A_{K}^{\#}$ of the birth-and-death process $(\NN_t)_{t\geq0}$ with transition rates \eqref{eq:TransitionRatesDual}.
\end{lemma}

\begin{proof}
	The result follows the same approach as in \cite{BarbourEthierGriffiths} (see p.\ 128). By Corollary 4.4.15 in \cite{EthierKurtz}, it suffices to show that there exists $H: \ZZ \mapsto [0, \infty)$ for which
	\begin{align}\label{eq:Lemma Condition1}
		h(\xx, \nn) + \left|\left(\A_Kh\right)\left(\xx, \nn\right)\right| &\leq H(\nn), & (\xx, \nn) \in \Delta_K \times \ZZ,
	\end{align}
	and
	\begin{align}\label{eq:LemmaCondition2}
		\mathcal{H} := \left\{ H\left(\NN_{t \wedge \tau_N}\right),\, 0 \leq t \leq t_0,\, N \geq 1 \right\} \textnormal{ is uniformly integrable}
	\end{align}
	for each initial state $\nn \in \ZZ$ and each $t_0 \geq 0$, where
\begin{equation}\nonumber
		\tau_N := \inf\left\{ s \geq 0 : \left|\NN_s\right| \geq N \right\}.
\end{equation} 
	Now, noting that $\binom{\left|\nn\right|}{\nn}\xx^{\nn} \leq 1$ and $e^{-\langle \ss, \xx\rangle} \leq e^{-\langle\min_{k}\sigma_k, \xx\rangle} = e^{-\min_{k}\sigma_k}$, we find that for $\xx\in\Delta_K, n \in \ZZ,$
	\begin{equation}\label{eq:LemmaDefnH}
		\begin{aligned}
			h(\xx, \nn) + \left|\left(\A_K h\right)\left(\xx, \nn\right)\right| &\leq \sum_{\lv\in\ZZ}\left(\delta_{\lv\nn} + \left|q(\nn, \lv)\right|\right)h(\xx, \lv) \\
			&\leq \sum_{\lv\in\ZZ}\left(\delta_{\lv\nn} + \left|q(\nn, \lv)\right|\right)\frac{e^{-\min_{k}\sigma_k}}{\binom{\left|\lv\right|}{\lv}\frac{B(\aa+\lv)}{c_{\aa, \ss}}}  
			=: H(\nn).
		\end{aligned}
	\end{equation} 
	We verify \eqref{eq:Lemma Condition1} and \eqref{eq:LemmaCondition2} for this choice of $H$. First note that
	\begin{equation}\nonumber
		\begin{aligned}
			\sum_{i=1}^{K}q(\nn, \nn - \boldsymbol{e}_i) &= \frac{\left|\nn\right|}{2}\left(\left|\aa\right| + \left|\nn\right|-1\right),  \\
			0 \leq \sum_{i=1}^{K}q(\nn, \nn + \boldsymbol{e}_i) &\leq \frac{\bar{\sigma}^{+}}{2}\left(\left|\aa\right| + \left|\nn\right|\right),
		\end{aligned}
	\end{equation} 
	where $\bar{\sigma}^{+} := \max_{k} \sigma_k = \max_{k} \left(\tilde{\sigma}_k - \min_{\ell}\tilde{\sigma}_{\ell}\right) = \max_k \tilde{\sigma}_k - \min_k \tilde{\sigma}_k$. Thus if $\varphi_{\ell} := \left(\left|\nn\right| + 1 \right)^{\ell}$ for $\ell \in \mathbb{Z}_{+}$, the generator $\A_K^{\#}$ of the dual birth-and-death process satisfies, as $\left|\nn\right| \rightarrow \infty$,
\begin{equation}\nonumber
\begin{aligned}
		\left(\A_K^{\#}\varphi_{\ell}\right)\left(\nn\right) \leq {}& \frac{\left|\nn\right|}{2}\left(\left|\nn\right| - 1 +\left|\aa\right|\right)\left[\left|\nn\right|^{\ell} - \left(\left|\nn\right| + 1\right)^{\ell}\right] \\
		{}& {}+ \frac{\bar{\sigma}^{+}}{2}\left(\left|\aa\right| + \left|\nn\right|\right)\left[\left(\left|\nn\right| + 2\right)^{\ell} - \left(\left|\nn\right| + 1\right)^{\ell}\right] \nonumber \\
		={} & {}-\frac{\ell}{2}\left|\nn\right|^{\ell+1} + \mathcal{O}\left(\left|\nn\right|^{\ell}\right)
		\leq C_{\ell},
\end{aligned}
\end{equation} 
	for some constant $C_{\ell}$ . Hence
	\begin{equation}\label{eq:LemmaExpectation}
		\begin{aligned}
			\mathbb{E}\left[\left(\left|\NN_{t \wedge \tau_N}\right| + 1\right)^{\ell} \big| \NN_0 = \nn\right] &= \left(\left|\nn\right| + 1\right)^{\ell} + \mathbb{E}\left[\int_{0}^{t\wedge \tau_N}\left(\A_K^{\#}\varphi_{\ell}\right)\left(\NN_{s}\right)ds\right]  \\
			&\leq \left(\left|\nn\right| + 1\right)^{\ell} + C_{\ell}t.
		\end{aligned}
	\end{equation} 
	Next for $\xx \in \Delta_K$ and $\nn \in \ZZ$,
\begin{equation}\nonumber
\begin{aligned}
		h(\xx, \nn) \leq \frac{e^{-\min_{k}\sigma_k}}{\binom{\left|\nn\right|}{\nn}\frac{B(\aa + \nn)}{c_{\aa, \ss}}} = \frac{e^{-\min_{k}\sigma_k}}{\binom{\left|\nn\right|}{\nn}}\frac{c_{\aa, \ss}\Gamma\left(\left|\nn\right| + \left|\aa\right|\right)}{\prod_{k=1}^{K}\Gamma\left(n_k + \sigma_k\right)} \leq B_0\left(\left|\nn\right| + 1\right)^{\ell_0}
\end{aligned}
\end{equation} 
	for constants $B_0 \geq 0$ and $\ell_0 \in \mathbb{Z}_{+}$ sufficiently large (using the same arguments as those employed in equation (3.22) of \cite{BarbourEthierGriffiths}). Hence for $H$ defined by \eqref{eq:LemmaDefnH} we have
\begin{align}\nonumber
		H(\nn) &\leq \sum_{\lv \in \ZZ}\left( \delta_{\lv\nn} + \left|q(\nn, \lv)\right| \right)e^{-\min_k \sigma_k}B_0\left(\left|\lv\right| + 1\right)^{\ell_0} \nonumber \\
		&\leq \sum_{\lv \in \ZZ}\left( \delta_{\lv\nn} + \left|q(\nn, \lv)\right| \right)B'_0\left(\left|\nn\right| + 2\right)^{\ell_0} \nonumber \\
		&\leq \sum_{\lv \in \ZZ}\left[ 1 + \left|\nn\right|\left(\left|\aa\right| + \left|\nn\right| - 1\right) + \bar{\sigma}^{+}\left(\left|\nn\right| + \left|\aa\right|\right) \right] B'_0\left(\left|\nn\right| + 2\right)^{\ell_0} \nonumber \\
		&\leq B''_0\left(\left|\nn\right|+2\right)^2\left(\left|\nn\right|+2\right)^{\ell_0} \nonumber \\
		&\leq B\left(\left|\nn\right| + 1\right)^{\ell_0 + 2}
\end{align}
	for sufficiently large constants $B'_0, B''_0, B \geq 0$. Thus by appealing to \eqref{eq:LemmaExpectation} the required uniform integrability condition \eqref{eq:LemmaCondition2} holds; in other words:
\begin{equation}\nonumber
\begin{aligned}
		\mathbb{E}\left[H\left(\NN_{t \wedge \tau_N}\right) \big| \NN_0 = \nn\right] \leq B\left[\left(\left|\nn\right| + 1\right)^{\ell_0 + 2} + C_{\ell_0 + 2}t_0\right] < \infty
\end{aligned}
\end{equation} 
	uniformly for $t \in [0, t_0]$ and $N \geq 1$, and if $\varepsilon > 0$ then there exists $\delta > 0$ such that, with an event $E$ satisfying $\mathbb{P}\left[E\right] < \delta$, we have
\begin{equation}\nonumber
\begin{aligned}
		\mathbb{E}\left[H\left(\NN_{t \wedge \tau_N}\right)1_{E} \big| \NN_0 = \nn\right] &\leq \mathbb{E}\left[B\left(\left|\NN_{t \wedge \tau_N}\right| + 1\right)^{\ell_0 + 2} 1_{E}\right] \nonumber \\
		&\leq B\left[\left(\left|\nn\right| + 1\right)^{\ell_0 + 2} + C_{\ell_0 + 2}t_0\right]\mathbb{P}\left[E\right] < \varepsilon.
\end{aligned}
\end{equation} 
	The latter inequality holds with the choice $\delta := \varepsilon \left[ B\left(\left|\nn\right| + 1\right)^{\ell_0 + 2} + C_{\ell_0 + 2}t_0 \right]^{-1}$.
\end{proof}

We refer to Proposition \ref{thm:Duality} as \emph{gamma duality}, in analogy with the terminology used in \cite{papas_ruggiero}, albeit in a completely different setting.
A key feature is that the dual rates in \eqref{eq:TransitionRatesDual} are fully explicit and do not involve ratios of intractable normalising constants. 
As a consequence, the dual process can be simulated efficiently (e.g.\ via standard Gillespie-type event simulation; see \cite{gillespie1977exact}) and its transition probabilities can be approximated or computed without repeatedly solving parameter-dependent integration problems.
This removes the main computational obstruction behind \eqref{eq:BEG_transition_expansion} and leads us to derive, in the next section, a Dirichlet-mixture representation for the transition function.


Next, we identify the stationary law of the $K$-type birth-and-death dual on $\ZZ$ with rates \eqref{eq:TransitionRatesDual}, focusing on the case $\aa \in (0,\infty)^K$.
Let
\[
\aa_{(\kk)} := \prod_{i=1}^K (\alpha_i)_{(k_i)},\qquad
\kk! := \prod_{i=1}^K k_i!,\qquad
|\aa|:=\sum_{i=1}^K \alpha_i,\qquad
|\aa|_{(m)} := (|\aa|)_{(m)},
\]
where $(\alpha)_{(k)} = \Gamma(\alpha + k)/\Gamma(\alpha)$ is the Pochhammer function.
\begin{proposition}\label{thm:StationaryDistDual}
	Let $(\NN_t)_{t\geq0}$ be as in Definition \ref{def:Dual} with $\aa \in (0,\infty)^K$. Then $(\NN_t)_{t\geq0}$ has reversible and stationary distribution
	\begin{equation}\label{eq:DualStationaryDistribution}
		\pi(\kk)
		\propto
		\frac{\aa_{(\kk)}}{|\aa|_{(|\kk|)}}\,
		\frac{\ss^{\kk}}{\kk!},
		\qquad \kk\in\ZZ,
	\end{equation}
	with normalising constant
	\begin{equation}\label{eq:normalisingConstant}
		\Phi^{(K)}_{2}(\aa;|\aa|;\ss)
		:=
		\sum_{\kk\in\ZZ}
		\frac{\aa_{(\kk)}}{|\aa|_{(|\kk|)}}\,
		\frac{\ss^{\kk}}{\kk!}.
	\end{equation} 
\end{proposition}

\begin{proof}
A simple calculation allows one to check that \eqref{eq:DualStationaryDistribution} satisfies the detailed balance condition, hence it is the unique stationary distribution for the  process in Definition \ref{def:Dual}. In particular, the fact that $\pi(\kk)q(\kk, \kk + \boldsymbol{e}_i)=\pi(\kk + \boldsymbol{e}_i)q(\kk + \boldsymbol{e}_i, \kk)$ follows easily by the identity
\begin{equation}\nonumber
\begin{aligned}
\frac{\aa_{(\kk)}}{\left(\left|\aa\right|\right)_{\left(\left|\kk\right|\right)}}\frac{\ss^{\kk}}{\kk!}
\sigma_i\left(\alpha_i + k_i\right) =
		\frac{\aa_{(\kk + \boldsymbol{e}_i)}}{\left(\left|\aa\right|\right)_{\left(\left|\kk + \boldsymbol{e}_i\right|\right)}}\frac{\ss^{\kk + \boldsymbol{e}_i}}{\left(\kk + \boldsymbol{e}_i\right)!}
		\left(k_i+1\right)\left(\left|\aa\right| + \left|\kk\right|\right),
\end{aligned}
\end{equation} 
while $\pi(\kk)q(\kk, \kk - \boldsymbol{e}_i)= \pi(\kk - \boldsymbol{e}_i)q(\kk - \boldsymbol{e}_i, \kk)$ follows from a similar computation, here omitted.
\end{proof}

The normalising constant $\Phi^{(K)}_{2}(\aa;|\aa|;\ss)$ is a  confluent Lauricella series (see, e.g., \citealp[Eq.\ (A.54), p.\ 215]{Mathai}), which is related to the normalising constant in \eqref{eq:StationaryWF} through
\begin{equation}
	\begin{aligned}
		\label{eq:c_as_Phi2}
		c_{\aa,\ss}
		=
		\int_{\Delta_K} e^{\langle \ss,\xx\rangle}\,\xx^{\aa-\mathbf 1}\,d\xx
		&=
		\sum_{\kk\in\ZZ}\frac{\ss^{\kk}}{\kk!}\int_{\Delta_K}\xx^{\aa+\kk-\mathbf 1}\,d\xx
		=
		\sum_{\kk\in\ZZ}\frac{\ss^{\kk}}{\kk!}B(\aa+\kk)
		\nonumber\\
		&=
		B(\aa)\sum_{\kk\in\ZZ}\frac{\aa_{(\kk)}}{|\aa|_{(|\kk|)}}\frac{\ss^{\kk}}{\kk!}
		=
		B(\aa)\,\Phi^{(K)}_{2}(\aa;|\aa|;\ss),
	\end{aligned}
\end{equation} 
where $B(\cdot)$ is the Beta function. Here we used $B(\aa+\kk)=B(\aa)\,\aa_{(\kk)}/|\aa|_{(|\kk|)}$ and Fubini's theorem to interchange integration and summation. We further note the following for later reference.

\begin{remark}\label{remark:DirichletConnection}
	If $\YY\sim\mathrm{Dir}(\aa)$ and $\tt\in[0,\infty)^K$, then
	\[
	\mathbb{E}\!\left[e^{\langle \tt,\YY\rangle}\right]
	=
	\Phi^{(K)}_{2}\!\left(\aa;|\aa|;\tt\right),
	\]
	so $\Phi^{(K)}_{2}(\aa;|\aa|;\cdot)$ is the mgf of a Dirichlet distribution with parameter $\aa$.
\end{remark}

Section~\ref{subsec:algorithmsStationaryDual} discusses how to efficiently evaluate $\Phi^{(K)}_{2}(\aa;|\aa|;\ss)$ and outlines a strategy to circumvent its evaluation altogether when simulating from \eqref{eq:DualStationaryDistribution}.


\subsection{Small-time behaviour}\label{subsec:smalltime}

A key ingredient in deriving a tractable transition function representation is that the dual process admits an entrance boundary at ``infinity'', which allows us to define an entrance law indexed by $\xx\in\Delta_K$. We first show that such entrance law can be constructed, and we subsequently characterise this law explicitly for the process in Definition \ref{def:Dual}.

\begin{proposition}\label{prop:ExtendToInfinity}
	Let $\rho : \ZZ \to \{0\}\cup(\mathbb{N}\times\Delta_K)$ be the one-to-one map defined by
\begin{equation}\nonumber
	\rho(\oo)=0,
		\qquad
		\rho(\nn)=\Big(|\nn|,\frac{\nn}{|\nn|}\Big)\quad \textnormal{if }\nn\neq\oo.
\end{equation} 
	Further, let $b_{\nn\lv}(t)=\mathbb{P}\left[\NN_t=\lv \mid \NN_0=\nn\right]$ denote the transition probabilities of the process $(\NN_t)_{t\ge0}$ as in Definition \ref{def:Dual}, and let $\ss \in [0,\infty)^K$. Then for $t>0$, $\lv\in\ZZ$ and $\xx\in\Delta_K$,
	\begin{align}\label{eq:StartFromInfinity}
		b_{\lv}(t;\xx)
		:=
		\liminf_{\rho(\nn)\to(\infty,\xx)} b_{\nn\lv}(t)
	\end{align}
	defines a probability mass function on $\ZZ$. Moreover, for any fixed $t>0$ and $\lv\in\ZZ$, the map $\yy\mapsto b_{\lv}(t;\yy)$ is weakly continuous on $\Delta_K$.
\end{proposition}

\begin{proof}
\cite{BarbourEthierGriffiths} provide very general conditions under which a continuous-time pure jump Markov process on $\mathbb{Z}_+^K$ can have an entrance boundary at infinity in the sense $\rho(\nn) \to (\infty,\xx)$. However, our setting does not satisfy their conditions, in particular due to the constant term in the birth rates in \eqref{eq:TransitionRatesDual}. Intuitively, the additional births slow down the entrance of the process from infinity, but this effect becomes negligible as $\left|\nn\right| \to \infty$ such that the desired conclusions still hold. In what follows we carefully verify that their strategy can be carried over to our case. Throughout the rest of the proof we will refer to their equations and sections as, e.g., BEG-(4.1) or BEG-Sec.~4.
	
In view of \eqref{eq:TransitionRatesDual}, $\NN_{t}$ has transitions 
$\nn \rightarrow \nn - \boldsymbol{e}_i $ at rate $n_i\left|\nn\right| + f_{1,i}(\nn) $ and $\nn \rightarrow \nn + \boldsymbol{e}_i $ at rate $\sigma_i\alpha_i + f_{2,i}(\nn) $, 
for $f_{1,i}(\nn) := 2q(\nn, \nn - \boldsymbol{e}_i) - n_i\left|\nn\right|$, and $f_{2, i}(\nn) := 2q(\nn, \nn+ \boldsymbol{e}_i) - \sigma_i\alpha_i$. It follows that $f_{1,i}(\oo) = 0 = f_{2,i}(\oo)$, 
\begin{equation}\nonumber
		f_{1, i}(\nn) = n_i(\left|\aa\right|-1) \geq -\left|\left|\aa\right|-1\right|n_i,
\end{equation} 
	and both functions are Lipschitz
\begin{equation}\nonumber
\begin{aligned}
		\left|f_{1,i}\left(\nn\right) - f_{1,i}\left(\nn + \boldsymbol{e}_i\right)\right| &= \left| n_i\left(\left|\aa\right|-1\right) - \left(n_i + 1\right)\left(\left|\aa\right|-1\right) \right| \leq \left|\left|\aa\right|-1\right|\left|\nn - \nn - \boldsymbol{e}_i\right|, \\
		\left|f_{2,i}\left(\nn\right) - f_{2,i}\left(\nn - \boldsymbol{e}_i\right)\right| &= \left| \sigma_in_i - \sigma_i(n_i - 1) \right| \leq \sigma_{i}\left|\nn - \left(\nn - \boldsymbol{e}_i\right)\right|,
\end{aligned}
\end{equation} 
	where we recall that $\sigma_i \geq 0$ for $i = 1, \dots, K$ in view of \eqref{eq:RescalingSigma_sec2}. 
Since we consider only genic selection, $\NN_{t}$ experiences only unit jumps at event times. 
	
	
As a consequence, the space-time generator $\mathscr{A}$ in BEG-(4.4) has an extra term (besides the absence of $f_{3,i}$), hence we have
\begin{equation}\nonumber
\begin{aligned}
		\left|\mathscr{A}U\left(\nn, u\right)\right| = &{} \left|2\left(u - \frac{1}{n}\right)\right| + \left|\left(n^2 + \sum_{i=1}^{K}f_{1,i}(\nn)\right)\left[\frac{1}{\left(n-1\right)^2} - \frac{1}{n^2} - 2u\left(\frac{1}{n-1} - \frac{1}{n}\right)\right]\right| \nonumber \\
		& {}+ \left|\left(S + \sum_{i=1}^{K}f_{2,i}(\nn)\right)\left[\frac{1}{\left(n+1\right)^2} - \frac{1}{n^2} - 2u\left(\frac{1}{n+1}  -\frac{1}{n}\right)\right]\right|  \nonumber \\
		\leq {}& 4\left(1 + c_1 + c_2 + S\right)\left(\frac{1}{n^2} + \frac{\left|u - \frac{1}{n}\right|}{n}\right),
\end{aligned}
\end{equation} 
	where $c_1 = K|\left|\aa\right| - 1|$, $c_2 = K\left|\ss\right|$, and $S := \langle\ss,\aa\rangle$. The bound on the extra term follows from
\begin{align}\nonumber
		S\bigg[&\frac{1}{(n+1)^2} - \frac{1}{n^2} -2u \bigg(\frac{1}{n+1} - \frac{1}{n}\bigg)\bigg] = S\left[\frac{-2n-1}{n^2(n+1)^2} + 2\left(u-\frac{1}{n}\right)\frac{1}{n(n+1)} + \frac{1}{n^2(n+1)}\right]  \\
		&= S\left[\frac{-1}{n(n+1)^2} + 2\left(u-\frac{1}{n}\right)\frac{1}{n(n+1)}\right]  \\
		&\leq S\left[\frac{1}{n(n+1)^2} + 2\left(u-\frac{1}{n}\right)\frac{2}{n}\right]  
		\leq 4S\left[\frac{1}{n^2} + \frac{\left|u - \frac{1}{n}\right|}{n}\right],
\end{align}
whereas the other bounds differ only in terms of the constants involved. Using these bounds, the results in BEG-Sec.~4 (chiefly Lemmas 4.2 and 4.5) hold in our case with different constants. We further point out that their Lemma 4.1 holds independently of the specific choice of space-time generator, and the proofs for Corollaries 4.4 and 4.6 can be applied directly to our case.
	
Similarly, BEG-(5.1) becomes
	\begin{align}\label{eq:BEG(5.1)}
		(\nn, \nn') &\rightarrow (\nn - \boldsymbol{e}_i, \nn' - \boldsymbol{e}_i) & \hphantom{}\nonumber \\
		&\hphantom{(\nn, \nn')} \textnormal{ at rate } \left(n_i \wedge n'_i\right)\left(\left|\nn\right| - C_1\right) + \tilde{f}_{1,i}(\nn) \wedge\tilde{f}_{1,i}(\nn') & \textnormal{for } & 1 \leq i \leq K \nonumber \\
		(\nn, \nn') &\rightarrow (\nn + \boldsymbol{e}_i, \nn' + \boldsymbol{e}_i) \textnormal{ at rate } \sigma_i\alpha_i + f_{2,i}(\nn)\wedge f_{2,i}(\nn') & \textnormal{for } & 1 \leq i \leq K \nonumber \\
		(\nn, \nn') &\rightarrow (\nn - \boldsymbol{e}_i, \nn' - \boldsymbol{e}_j) \textnormal{ at rate } \pi_{ij}(\nn, \nn')\left(\left|\nn\right| - C_1\right) & \textnormal{for } & 1 \leq i \neq j \leq K, 
	\end{align}
	where $C_1 := \left|\left|\aa\right| - 1\right|$, $\tilde{f}_{1,i}(\nn) := f_{1,i}(\nn) + C_1 n_i \geq 0$, and the $\pi_{ij}$ are chosen such that
\begin{equation}\nonumber
\begin{aligned}
		\sum_{j\neq i}\pi_{ij}\left(\nn, \nn'\right) = n_i - \left(n_i \wedge n_i'\right) \textnormal{ for } 1 \leq i \leq K \\
		\sum_{i\neq j}\pi_{ij}\left(\nn, \nn'\right) = n_j' - \left(n_j \wedge n_j'\right) \textnormal{ for } 1 \leq j \leq K.
\end{aligned}
\end{equation} 
Furthermore, the transitions relevant in our setting and corresponding to BEG-(5.13) become
\begin{equation}\nonumber
\begin{aligned}
		(\nn, \nn') &\rightarrow (\nn - \boldsymbol{e}_i, \nn') & \textnormal{ at rate } & \tilde{f}_{1,i}(\nn) - \tilde{f}_{1,i}(\nn) \wedge\tilde{f}_{1,i}(\nn') & \textnormal{for } & 1 \leq i \leq K, \nonumber\\
		(\nn, \nn') &\rightarrow (\nn, \nn' - \boldsymbol{e}_i) & \textnormal{ at rate } & \tilde{f}_{1,i}(\nn') - \tilde{f}_{1,i}(\nn) \wedge\tilde{f}_{1,i}(\nn') & \textnormal{for } & 1 \leq i \leq K,\nonumber\\
		(\nn, \nn') &\rightarrow (\nn + \boldsymbol{e}_i, \nn') & \textnormal{ at rate } & \tilde{f}_{2,i}(\nn) - \tilde{f}_{2,i}(\nn) \wedge\tilde{f}_{2,i}(\nn') & \textnormal{for } & 1 \leq i \leq K,\nonumber\\
		(\nn, \nn') &\rightarrow (\nn, \nn' + \boldsymbol{e}_i) & \textnormal{ at rate } & \tilde{f}_{2,i}(\nn) - \tilde{f}_{2,i}(\nn) \wedge\tilde{f}_{2,i}(\nn') & \textnormal{for } & 1 \leq i \leq K.
\end{aligned}
\end{equation} 
	Here we again need to consider the appropriate terms appearing in the space-time generator $\mathscr{A}$ to apply the same proofs for Lemmas 5.1 and 5.3, thereby leading to the same qualitative result albeit with different constants. Lemma 5.2, Corollary 5.4, and Theorems 5.5 and 5.6 apply directly to our case from which we deduce that \eqref{eq:StartFromInfinity} defines a probability distribution on $\ZZ$ as desired. Moreover, Theorem 5.7  also applies directly to give us the weak continuity of $b_{\lv}(t; \yy)$ in $\yy \in \Delta_K$.
\end{proof}


Let $(\NN_t)_{t\ge0}$ be as in Definition \ref{def:Dual}, and write
\[
N_t := |\NN_t|=\sum_{i=1}^K N_{t,i}.
\]
Then $N_t$ has jumps of size $\pm 1$, with death rate
\begin{align}\nonumber
	\mu_n := \frac{n(n+|\aa|-1)}{2},
\end{align}
while the birth rate depends on the current composition: if $\NN_t=\nn$ with $|\nn|=n$, then
\begin{equation}\nonumber
	\lambda(\nn) := \frac{\langle \ss,\aa + \nn\rangle}{2},
\end{equation} 
To describe the entrance regime relevant for \eqref{eq:StartFromInfinity}, fix $\xx\in\Delta_K$ and set
$\tilde s:=\langle \ss,\xx\rangle$.
Under the entrance scaling $\rho(\nn)\to(\infty,\xx)$, by definition $n:=|\nn|\to\infty$ and $\nn/n\to\xx$.
Therefore
\[
\langle \ss,\nn\rangle
=
n\Big\langle \ss,\frac{\nn}{n}\Big\rangle
=
n\langle \ss,\xx\rangle
+
n\Big\langle \ss,\frac{\nn}{n}-\xx\Big\rangle
=
n\tilde s+o(n),
\]
and hence $\lambda(\nn)$ can be written, for $S = \langle\ss, \aa\rangle$,
\[
\lambda_{n}
=
\frac{S+n\tilde s+o(n)}{2},
\qquad n\to\infty.
\]

The following theorem extends \citet[Theorem 4]{Griffiths84}.
\begin{theorem}\label{thm:CLT}
	Let $(M_t)_{t\ge0}$ be a birth-and-death process with rates $\lambda_{n},\mu_{n}$ as above. 
	Define
	\begin{align}\label{eq:v(t)}
		v(t):=\frac{r_{+}-r_{-}e^{-(r_{+}-r_{-})t/2}}{1-e^{-(r_{+}-r_{-})t/2}},
		\qquad
		r_{\pm}:=\frac{(\tilde s-|\aa|+1)\pm\sqrt{(\tilde s-|\aa|+1)^2+4S}}{2}.
	\end{align}
	Then, as $t\downarrow  0$,
	\begin{align}\label{eq:CLT}
		\sqrt{\frac{3}{v(t)}}\bigl(M_t-v(t)\bigr)\ \overset{d}{\longrightarrow}\ Z,
	\end{align}
	with $Z\sim \mathcal{N}(0,1)$.
\end{theorem}

\begin{proof}
Similarly to the proof of Proposition \ref{prop:ExtendToInfinity}, existing theory does not seem to cover our case due to the constant terms in the birth rates in \eqref{eq:TransitionRatesDual}, which in turn determine the absence of an absorbing state and the existence of a non-degenerate limiting distribution on $\mathbb{Z}_{+}$. Nonetheless, the small-time dynamics are still dominated by the quadratic death rates, and hence we adapt the arguments used in \citet[Theorem 5.1]{Bansaye}, cited henceforth as BMR. 

Observe that $\lim_{n\rightarrow\infty}\frac{\lambda_n}{\mu_n} = \lim_{n\rightarrow\infty} \frac{S+\tilde{s}n}{n(n+|\aa|-1)} = 0$, and furthermore
\begin{equation}\nonumber
		\sum_{n\geq1}\frac{1}{n}\frac{\lambda_n}{\mu_n} = \sum_{n\geq1}\frac{n\tilde{s} + S}{n^{2}\left(n + |\aa| - 1\right)} < \infty.
\end{equation} 
	Recall that a non-zero sequence of numbers $(x_n)_{n\geq0}$ varies regularly with index $\rho \neq 0$ if for any $a > 0$
\begin{equation}\nonumber
		\lim_{n\rightarrow\infty}\frac{x_{[an]}}{x_n} = a^{\rho}.
\end{equation} 
	From this it follows immediately that $(\mu_n)_{n\geq0}$ varies regularly with index $2$, and thus by BMR-Theorem 5.1 we can conclude that 
	\begin{equation}\nonumber
\begin{aligned}
		\sqrt{\left(2\rho - 1\right)v(t)}\left(\frac{M_t}{v(t)} - 1\right) \overset{t\rightarrow0}{\longrightarrow} Z
\end{aligned}
\end{equation} 
	where $Z \sim \mathcal{N}(0,1)$, and $v(t)$ is regularly varying with index $(1-\rho)^{-1}$ for some index $\rho > 1$.

To characterise the regularly varying function $v(t)$ explicitly, denote now by $\mathbb{P}_n$ the law of the process $(M_t)_{t\geq0}$ started from $n\in\mathbb{N}$ and define $\mathbb{E}_{n}$ to be the expectation with respect to this law. By BMR-Lemma 2.1, $(\mathbb{P}_n)_{n\in\mathbb{N}}$ has a limit $\mathbb{P}$, and we define $\mathbb{E}$ to be the expectation taken with respect to $\mathbb{P}$. As defined in BMR-Section 4 we have that
\begin{equation}\nonumber
		v(t) := \inf\left\{ n \geq 0 : \mathbb{E}\left[T_n\right] \leq t \right\}
\end{equation} 
	where $T_n = \inf\{t\geq 0:\: M_t = n\}$. From this we deduce that $v(t)$ for small $t$ should solve the ODE
	\begin{align}\label{eq:vtODE}
		v'(t) &= -\frac{v(t)\left(v(t) + |\aa| - 1\right)}{2} + \frac{S + \tilde{s}v(t)}{2}.
	\end{align}
	Briefly, define $G(n) := \sum_{k\geq n}\mathbb{E}_{k+1}\left[T_k\right]$, and observe that $G(n)$ is deterministic and strictly non-increasing. Further, we see that $v(t)$ is the inverse mapping of $G(n)$ and by definition $G(v(t)) \leq t < G(v(t) - 1)$, such that $G(v(t)) \approx t$. Additionally, 
\begin{equation}\nonumber
		\mathbb{E}_{n+1}\left[T_n\right] = \frac{1}{\lambda_n + \mu_n} + \frac{\lambda_n}{\lambda_n + \mu_n}\left(\mathbb{E}_{n+2}\left[T_{n+1}\right] + \mathbb{E}_{n+1}\left[T_{n}\right]\right) 
\end{equation} 
	from which it follows that $\mathbb{E}_{n+1}\left[T_n\right] = \frac{1}{\mu_n} + \frac{\lambda_n}{\mu_n}\mathbb{E}_{n+2}\left[T_{n+1}\right]$. Since we are considering small time $t\rightarrow0$ and thus $n\rightarrow\infty$, we have that $\mathbb{E}_{n+2}\left[T_{n+1}\right] \approx \mathbb{E}_{n+1}\left[T_{n}\right]$, and thus $\mathbb{E}_{n+1}\left[T_n\right] \approx (\mu_n - \lambda_n)^{-1}$. Combining this with $G'(v(t))v'(t) \approx 1$ leads to the ODE above.
	
	Solving \eqref{eq:vtODE} gives \eqref{eq:v(t)}.
	Now observe that
\begin{equation}\nonumber
		\frac{v(at)}{v(t)} = \frac{r_{-} + \frac{r_{+} - r_{-}}{1-e^{-\frac{\left(r_{+} - r_{-}\right)at}{2}}}}{r_{-} + \frac{r_{+} - r_{-}}{1-e^{-\frac{\left(r_{+} - r_{-}\right)t}{2}}}} \sim \frac{1-e^{-\frac{\left(r_{+} - r_{-}\right)t}{2}}}{1-e^{-\frac{\left(r_{+} - r_{-}\right)at}{2}}},
\end{equation} 
	as $t\rightarrow 0$. Let $g(x) := \frac{1 - e^{-x}}{x}$, then $g(x) \rightarrow 1$ as $x \rightarrow 0$ such that
\begin{equation}\nonumber
\begin{aligned}
		\frac{1-e^{-\frac{\left(r_{+} - r_{-}\right)t}{2}}}{1-e^{-\frac{\left(r_{+} - r_{-}\right)at}{2}}} = \frac{g\left(\left(\frac{r_{+}  -r_{-}}{2}\right)t\right)}{ag\left(\left(\frac{r_{+}  -r_{-}}{2}\right)at\right)} \overset{t\rightarrow0}{\longrightarrow} a^{-1},
\end{aligned}
\end{equation} 
	and thus we conclude that $v(t)$ is indeed regularly varying at 0 with index $-1$.
\end{proof}

The centering $v(t)$ captures the explosion scale of the entrance regime. As $t\downarrow 0$,
\[
v(t)=\frac{2}{t}+\frac{r_{+}+r_{-}}{2}
+\Big(\frac{r_{+}-r_{-}}{2}\Big)^2\frac{t}{6}+o(t).
\]
Hence the leading term $2/t$ agrees with the classical coming-down-from-infinity behaviour for the block-counting process of Kingman's coalescent with mutation \citep{Griffiths84,Bansaye}, while the effect of selection enters only through the lower-order corrections.

From a computational perspective, Theorem~\ref{thm:CLT} provides a sharp description of the growth of the total count $N_t=|\NN_t|$ in the entrance regime. In particular, the approximation
$N_t\approx v(t)+\sqrt{v(t)/3}\,Z$ for small $t$ yields a principled approximation for simulating the dual from its entrance boundary by first drawing a plausible $N_t$ and then sampling the composition conditional on $N_t$. We exploit these approximations in Section~\ref{subsec:algorithmsWFDiff}.


\section{A new transition function representation}\label{sec:transfunc}


\subsection{A tractable expansion with Dirichlet components}

Combining Proposition \ref{thm:Duality} with Proposition~\ref{prop:ExtendToInfinity}
and a sample-inversion argument (see \cite{BarbourEthierGriffiths,EG09,GriffithsJenkinsLessard,griffiths2025dual})
yields the following Dirichlet-mixture representation of the non-neutral transition function in terms of the dual started from its entrance boundary.

\begin{theorem}\label{thm:TransitionDensityWF}
	The transition function of the $K$-dimensional WF diffusion with generator \eqref{eq:WFGeneratorHaploidPIM} admits the representation
	\begin{align}\label{eq:NonNeutralTransitionDensityWF}
		P_{\aa}^{\ss}(\xx,d\yy;t)
		=
		\sum_{\lv\in\ZZ}
		b_{\lv}(t;\xx)\,
		\Pi^{\oo}_{\aa+\lv}(d\yy),
	\end{align}
	where $(\NN_t)_{t\ge0}$ is as given in Definition \ref{def:Dual}, $\ss$ is as in \eqref{eq:RescalingSigma_sec2}, and $b_{\lv}(t;\xx)$ is defined in \eqref{eq:StartFromInfinity}.
\end{theorem}

\begin{proof}
	First we assume $\aa \in (0,\infty)^K$. Recall that from Proposition \ref{thm:Duality} we have that given $\xx \in \Delta_K, \mm \in \ZZ$ and any $t > 0$
\begin{equation}\nonumber
		\mathbb{E}\left[h\left(\XX_t, \nn\right) \big| \XX_0 = \xx\right] = \mathbb{E}\left[h\left(\xx, \NN_t\right) \big| \NN_0 = \nn\right],
\end{equation} 
	where $h(\xx, \nn) = \frac{c_{\aa, \ss}}{c_{\aa + \nn, \oo}}e^{-\langle\ss, \xx\rangle}\xx^{\nn}$. Since we can recover a distribution from its moments (see equation (2.13) and subsequent argument in \cite{BarbourEthierGriffiths}, but also \cite{EG09,GriffithsJenkinsLessard}), we want to adapt a similar approach whilst explicitly accounting for the fact that $h(\xx, \nn)$ gives access only to exponentially tilted moments. To this end, let $\mathscr{P}(\Delta_K)$ denote the space of all Borel probability measures on $\Delta_K$, and for $\mu \in \mathscr{P}\left(\Delta_K\right)$, define $\tilde{\Lambda}_n : \mathscr{P}\left(\Delta_K\right) \rightarrow \mathscr{P}\left(\Delta_K\right)$ by
\begin{equation}\nonumber
\begin{aligned}
		\tilde{\Lambda}_n \mu := \sum_{\mm\in\ZZ: \left|\mm\right| = n}\binom{\left|\mm\right|}{\mm}\left[\int_{\Delta_K}\tilde{h}(\boldsymbol{z}, \mm)\mu(d\boldsymbol{z})\right]\delta_{\mm / \left|\mm\right|}.
\end{aligned}
\end{equation} 
	For $G \subseteq \Delta_K$ open we have
\begin{equation}\label{eq:ThmTildeLambdaNG}
\begin{aligned}
		\liminf\limits_{n\rightarrow \infty}\tilde{\Lambda}_n \mu\left(G\right) &= \liminf_{n\rightarrow\infty}\int_{\Delta_K}\sum_{\mm\in\ZZ : \left|\mm\right| = n}\binom{\left|\mm\right|}{\mm}\tilde{h}(\boldsymbol{z}, \mm)\delta_{\mm / \left|\mm\right|}\left(G\right)\mu(d\boldsymbol{z})  \\
		&\geq \int_{\Delta_K}\liminf_{n\rightarrow\infty}\sum_{\mm\in\ZZ : \left|\mm\right| = n}\binom{\left|\mm\right|}{\mm}\tilde{h}(\boldsymbol{z}, \mm)\delta_{\mm / \left|\mm\right|}\left(G\right)\mu(d\boldsymbol{z})  \\
		&= \int_{\Delta_K}\liminf_{n\rightarrow\infty}\sum_{\mm\in\ZZ : \left|\mm\right| = n}\binom{\left|\mm\right|}{\mm}\boldsymbol{z}^{\mm}\delta_{\mm / \left|\mm\right|}\left(G\right)\tilde{\mu}(d\boldsymbol{z}) = \tilde{\mu}\left(G\right)
\end{aligned}
\end{equation} 
	where $\tilde{\mu}$ is defined via
\begin{equation}\nonumber
		\frac{d\tilde{\mu}}{d\mu}\left(\yy\right) = e^{-\langle\ss, \yy\rangle}.
\end{equation} 
	Equation \eqref{eq:ThmTildeLambdaNG} implies that $\tilde{\Lambda}_n\mu \overset{\textnormal{d}}{\rightarrow} \tilde{\mu}$ as $n \rightarrow \infty$ by the Portmanteau theorem, thus we do not recover the distribution of interest $\mu$ but rather an exponentially tilted version. In view of Proposition \ref{prop:ExtendToInfinity}, we have that for any $\xx \in \Delta_K$ and $t > 0$
\begin{equation}\nonumber
\begin{aligned}
		b_{\lv}(t; \xx) &:= \liminf\limits_{\rho(\mm) \rightarrow \left(\infty, \xx\right)}b_{\mm\lv}(t)
\end{aligned}
\end{equation} 
	defines a probability distribution on $\ZZ$. 
	
	Now fix $\xx \in \Delta_K$ and $t > 0$, then using \eqref{eq:GammaDuality}
	\begin{equation}\nonumber
		\begin{aligned}
			\mathbb{E}\left[h\left(\XX_t, \mm\right) \big| \XX_0 = \xx\right] &= \frac{\mathbb{E}\left[\tilde{h}\left(\XX_t, \mm\right) \big| \XX_0 = \xx\right]}{\mathbb{E}\left[\tilde{h}\left(\XX_\infty, \mm\right) \right]}	 \\
			&= \mathbb{E}\left[\frac{\tilde{h}(\xx, \NN_t)}{\mathbb{E}\left[\tilde{h}(\XX_{\infty}, \NN_t) \big| \NN_t\right]}\big| \NN_0 = \mm\right] 
			= \mathbb{E}\left[h\left(\xx, \NN_t\right) \big| \NN_0 = \mm\right],
		\end{aligned}
	\end{equation} 
	which gives us that
	\begin{equation}\label{eq:ThmPhiDefn}
		\begin{aligned}
			\liminf\limits_{\rho(\mm) \rightarrow \left(\infty, \yy\right)}\frac{\mathbb{E}\left[\tilde{h}\left(\XX_t, \mm\right) \big| \XX_0 = \xx\right]}{\mathbb{E}\left[\tilde{h}\left(\XX_\infty, \mm\right) \right]} &= \liminf\limits_{\rho(\mm) \rightarrow \left(\infty, \yy\right)} \mathbb{E}\left[\frac{\tilde{h}(\xx, \NN_t)}{\mathbb{E}\left[\tilde{h}(\XX_{\infty}, \NN_t) \big| \NN_t\right]}\big| \NN_0 = \mm\right]  \\
			&=\liminf\limits_{\rho(\mm) \rightarrow \left(\infty, \yy\right)} \sum_{\lv\in\ZZ} b_{\mm\lv}(t)\frac{\tilde{h}(\xx, \lv)}{\mathbb{E}\left[\tilde{h}(\XX_\infty, \lv)\right]}  \\
			&\geq \sum_{\lv\in\ZZ}b_{\lv}(t; \yy)\frac{\tilde{h}(\xx, \lv)}{\mathbb{E}\left[\tilde{h}(\XX_\infty, \lv)\right]} 
			=: \phi(t, \xx, \yy).
		\end{aligned}
	\end{equation} 
	Now define
	\begin{align}\label{eq:ThmTildePiP}
		\tilde{\Pi}(d\yy) &:= e^{-\langle\ss, \yy\rangle}\Pi_{\aa}^{\ss}(d\yy), & \tilde{P}_{t, \xx}(d\yy) &:= e^{-\langle\ss, \yy\rangle}P_{\aa}^{\ss}(\xx, d\yy; t), \nonumber \\
		\tilde{\mathbb{P}}_0(d\yy) &:= \frac{\tilde{\Pi}(d\yy)}{\tilde{\Pi}(\Delta_K)}, & \tilde{\mathbb{P}}_{t, \xx}(d\yy) &:= \frac{\tilde{P}_{t, \xx}(d\yy)}{\tilde{P}_{t, \xx}(\Delta_K)},
	\end{align}
	and since
\begin{equation}\nonumber
\begin{aligned}
		\frac{\mathbb{E}\left[\tilde{h}\left(\XX_t, \mm\right) \big| \XX_0 = \xx\right]}{\mathbb{E}\left[\tilde{h}\left(\XX_\infty, \mm\right) \right]} = \frac{\tilde{P}_{t, \xx}(\Delta_K)}{\tilde{\Pi}(\Delta_K)}\frac{\int_{\Delta_K}\boldsymbol{z}^{\mm}\tilde{\mathbb{P}}_{t, \xx}(d\boldsymbol{z})}{\int_{\Delta_K}\boldsymbol{z}^{\mm}\tilde{\mathbb{P}}_0(d\boldsymbol{z})},
\end{aligned}
\end{equation} 
	we get that \eqref{eq:ThmPhiDefn} becomes
\begin{equation}\nonumber
\begin{aligned}
		\liminf\limits_{\rho(\mm) \rightarrow \left(\infty, \yy\right)} \frac{\int_{\Delta_K}\boldsymbol{z}^{\mm}\tilde{\mathbb{P}}_{t, \xx}(d\boldsymbol{z})}{\int_{\Delta_K}\boldsymbol{z}^{\mm}\tilde{\mathbb{P}}_0(d\boldsymbol{z})} \geq \frac{\tilde{\Pi}(\Delta_K)}{\tilde{P}_{t, \xx}(\Delta_K)} \phi(t, \xx, \yy).
\end{aligned}
\end{equation} 
	Assume now that
	\begin{align}\label{eq:ThmLemma2.1Condition}
		\frac{\tilde{\Pi}(\Delta_K)}{\tilde{P}_{t, \xx}(\Delta_K)} \int_{\Delta_K}\phi(t, \xx, \yy)\tilde{\mathbb{P}}_{0}(d\yy) = 1,
	\end{align}
	then we can apply Lemma 2.1 in \cite{BarbourEthierGriffiths} with $\mu = \tilde{P}_{t,\xx}$, $\mu_0 = \tilde{\mathbb{P}}_0$ and $\varphi(\yy) = \frac{\tilde{\Pi}(\Delta_K)}{\tilde{P}_{t, \xx}(\Delta_K)} \phi(t, \xx, \yy)$ to get that
\begin{equation}\nonumber
\begin{aligned}
		\tilde{\mathbb{P}}_{t, \xx}(d\yy) = \frac{\tilde{\Pi}(\Delta_K)}{\tilde{P}_{t, \xx}(\Delta_K)} \phi(t, \xx, \yy)\tilde{\mathbb{P}}_{0}(d\yy)
\end{aligned}
\end{equation} 
	for $\Pi^{\ss}_{\aa}$-a.e.\ $\xx \in \Delta_K$. Re-arranging the above, we get that
\begin{equation}\nonumber
		P_{\aa}^{\ss}(\xx, d\yy; t) = \phi(t, \xx, \yy)\Pi^{\ss}_{\aa}(d\yy)
\end{equation} 
	for $\Pi^{\ss}_{\aa}$-a.e.\ $\xx \in \Delta_K$, and since $P_{\aa}^{\ss}(\xx, d\yy; t)$ is $\Pi^{\ss}_{\aa}$-reversible, 
\begin{equation}\nonumber
		\Pi^{\ss}_{\aa}(d\xx)P_{\aa}^{\ss}(\xx, d\yy; t) = \Pi^{\ss}_{\aa}(d\yy)P_{\aa}^{\ss}(\yy, d\xx; t),
\end{equation} 
	we conclude that for $\Pi^{\ss}_{\aa} \times \Pi^{\ss}_{\aa}$-a.e.\ $(\xx, \yy) \in \Delta_K \times \Delta_K$ we have $\phi(t, \xx, \yy) = \phi(t, \yy, \xx)$ and hence
\begin{equation}\nonumber
\begin{aligned}
		P_{\aa}^{\ss}(\xx, d\yy; t) 
		&= \phi(t, \yy, \xx)\Pi^{\ss}_{\aa}(d\yy) 
		= \sum_{\lv\in\ZZ}b_{\lv}(t; \xx)\frac{h(\yy, \lv)}{\mathbb{E}\left[h(\XX_{\infty}, \lv)\right]}\Pi^{\ss}_{\aa}(d\yy) \nonumber \\
		&= \sum_{\lv\in\ZZ}b_{\lv}(t; \xx)\frac{\yy^{\lv}e^{-\langle\ss, \yy\rangle}}{B(\aa + \lv) / c_{\aa, \ss}}\frac{\yy^{\aa - \boldsymbol{1}}e^{\langle\ss, \yy\rangle}}{c_{\aa, \ss}}d\yy \nonumber \\
		&= \sum_{\lv\in\ZZ}b_{\lv}(t; \xx)\frac{\yy^{\aa + \lv - \boldsymbol{1}}}{B(\aa + \lv)}d\yy 
		= \sum_{\lv\in\ZZ}b_{\lv}(t; \xx)\Pi^{\oo}_{\aa+\lv}(d\yy),
\end{aligned}
\end{equation} 
	as required. It thus remains to prove that \eqref{eq:ThmLemma2.1Condition} holds.
	
	To this end, define $\tilde{h}_n : \Delta_K \rightarrow [0, \infty)$ and $\tilde{\phi}_n : \Delta_K \rightarrow [0, \infty)$ by
\begin{equation}\nonumber
\begin{aligned}
		\tilde{h}_n(\yy) &:= \begin{cases}
			\frac{\int_{\Delta_K}\tilde{h}(\boldsymbol{z}, \mm)P_{\aa}^{\ss}(\xx, d\boldsymbol{z}; t)}{\int_{\Delta_K} \tilde{h}(\boldsymbol{z}, \mm)\Pi^{\ss}_{\aa}(d\boldsymbol{z})} & \textnormal{ if } \yy = \mm/|\mm|, \textnormal{ where } \mm \in \ZZ \textnormal{ and } \left|\mm\right| = n \\
			0 & \textnormal{otherwise},
		\end{cases} \nonumber \\
		\tilde{\phi}_n(\yy) &:= \sum_{\mm\in \ZZ : \left|\mm\right| = n} \tilde{h}_n\left(\frac{\mm}{\left|\mm\right|}\right)\binom{\left|\mm\right|}{\mm}\tilde{h}(\yy, \mm) = \int_{\Delta_K}\tilde{h}_n(\boldsymbol{z})\tilde{\Lambda}_n\delta_{\yy}(d\boldsymbol{z}).
\end{aligned}
\end{equation} 
	Now note that for each $n \in \mathbb{Z}_{+}$
\begin{equation}\label{eq:ThmLemma2.1PartI}
\begin{aligned}
		\int_{\Delta_K} \tilde{\phi}_n(\yy)e^{\langle \ss, \yy \rangle}\tilde{\mathbb{P}}_0(d\yy) &= \frac{1}{\tilde{\Pi}(\Delta_K)}\int_{\Delta_K}\tilde{\phi}_n(\yy)\Pi^{\ss}_{\aa}(d\yy)  \\
		&= \frac{1}{\tilde{\Pi}(\Delta_K)}\int_{\Delta_K}\sum_{\mm\in \ZZ : \left|\mm\right| = n} \tilde{h}_n\left(\frac{\mm}{\left|\mm\right|}\right)\binom{\left|\mm\right|}{\mm}\tilde{h}(\yy, \mm)\Pi^{\ss}_{\aa}(d\yy)  \\
		&= \frac{1}{\tilde{\Pi}(\Delta_K)}\sum_{\mm\in \ZZ : \left|\mm\right| = n} \binom{\left|\mm\right|}{\mm}\int_{\Delta_K}\tilde{h}(\yy, \mm)P_{\aa}^{\ss}(\xx, d\yy; t)  \\
		&= \frac{\tilde{P}_{t, \xx}(\Delta_K)}{\tilde{\Pi}(\Delta_K)}\sum_{\mm\in \ZZ : \left|\mm\right| = n} \binom{\left|\mm\right|}{\mm}\int_{\Delta_K}\yy^{\mm}\tilde{P}_{t, \xx}(d\yy) 
		= \frac{\tilde{P}_{t, \xx}(\Delta_K)}{\tilde{\Pi}(\Delta_K)},
\end{aligned}
\end{equation} 
	and for $\yy \in \Delta_K$
	\begin{align}\label{eq:ThmLemma2.1PartII}
		\liminf_{n\rightarrow\infty}\tilde{\phi}_n(\yy) = \liminf_{n\rightarrow\infty}\int_{\Delta_K}\tilde{h}_n(\boldsymbol{z})\tilde{\Lambda}_n\delta_{\yy}(d\boldsymbol{z}) \geq \int_{\Delta_K}\phi(t, \xx, \boldsymbol{z})\tilde{\delta}_{\yy}(d\boldsymbol{z}) = \phi(t, \xx, \yy)e^{-\langle\ss, \yy\rangle},
	\end{align}
	where we make use of Fatou's lemma for the inequality, \eqref{eq:ThmPhiDefn} for the $\liminf$ of $\tilde{h}_n$, \eqref{eq:ThmTildeLambdaNG} for the $\liminf$ of $\tilde{\Lambda}_n\delta_{\yy}$, and the Skorokhod representation theorem to put everything onto a single probability space. Here $\tilde{\delta}_{\yy}$ is the measure defined via 
\begin{equation}\nonumber
		\frac{d\tilde{\delta}_{\yy}}{d\delta_{\yy}}(\boldsymbol{z}) = e^{-\langle \ss, \boldsymbol{z} \rangle}
\end{equation} 
	for $\boldsymbol{z} \in \Delta_K$. 
	
	Hence combining \eqref{eq:ThmLemma2.1PartI} with \eqref{eq:ThmLemma2.1PartII} and again using Fatou's lemma allows us to deduce that for any $t > 0$ and $\xx \in \Delta_K$
	\begin{align}
		\frac{\tilde{P}_{t, \xx}(\Delta_K)}{\tilde{\Pi}(\Delta_K)} &= \liminf_{n\rightarrow\infty}\int_{\Delta_K}\tilde{\phi}_n(\yy)e^{\langle\ss, \yy\rangle}\tilde{\mathbb{P}}_{0}(d\yy) \nonumber \\
		&\geq \int_{\Delta_K}\liminf_{n\rightarrow\infty}\tilde{\phi}_n(\yy)e^{\langle\ss, \yy\rangle}\tilde{\mathbb{P}}_{0}(d\yy) \geq \int_{\Delta_K}\phi(t, \xx, \yy)\tilde{\mathbb{P}}_0(d\yy).	\label{eq:ThmLemma2.1PartIII}
	\end{align}
	If we can show that the above inequality is in fact an equality then we are done. To this end, consider the probability measure $\mathbb{Q}(d\xx)$ defined by
\begin{equation}\nonumber
		\mathbb{Q}(d\xx) = \frac{\tilde{P}_{t, \xx}(\Delta_K)\Pi^{\ss}_{\aa}(d\xx)}{\int_{\Delta_K}\tilde{P}_{t, \boldsymbol{z}}(\Delta_K)\Pi^{\ss}_{\aa}(d\boldsymbol{z})} = \frac{\tilde{P}_{t, \xx}(\Delta_K)}{\tilde{\Pi}(\Delta_K)}\Pi^{\ss}_{\aa}(d\xx),
\end{equation} 
	where the expression for the normalising constant follows from the fact that $\Pi^{\ss}_{\aa}$ is the stationary measure for $\XX$. 
	
	Now for fixed $t$ we have that
\begin{equation}\nonumber
\begin{aligned}
		\int_{\Delta_K}&\int_{\Delta_K}\phi(t, \xx, \yy)\frac{\tilde{\Pi}(\Delta_K)}{\tilde{P}_{t, \xx}(\Delta_K)}\tilde{\mathbb{P}}_0(d\yy)\mathbb{Q}(d\xx) 
		= \int_{\Delta_K}\int_{\Delta_K}\phi(t, \xx, \yy)\frac{\tilde{\Pi}(\Delta_K)}{\tilde{P}_{t, \xx}(\Delta_K)}\mathbb{Q}(d\xx)\tilde{\mathbb{P}}_0(d\yy) \nonumber \\
		&= \int_{\Delta_K}\int_{\Delta_K}\phi(t, \xx, \yy)\Pi^{\ss}_{\aa}(d\xx)\tilde{\mathbb{P}}_0(d\yy) \nonumber 
		= \int_{\Delta_K}\sum_{\lv\in\ZZ}b_{\lv}(t; \yy)\tilde{\mathbb{P}}_0(d\yy) = 1
\end{aligned}
\end{equation} 
	where we use the fact that $\{b_{\lv}(t; \xx) : \lv \in \ZZ\}$ defines a probability distribution on $\ZZ$ as established in Proposition \ref{prop:ExtendToInfinity}. This implies that \eqref{eq:ThmLemma2.1PartIII} is indeed an equality, at least for $\mathbb{Q}$-a.e.\ $\xx\in\Delta_K$. But since $\XX$ is Feller, $\tilde{P}_{t, \xx}(\Delta_K)$ is continuous as a function of $\xx$, so $\mathbb{Q} \sim \Pi^{\ss}_{\aa}$ and we can replace ``$\mathbb{Q}$-a.e.'' with ``$\Pi^{\ss}_{\aa}$-a.e.''. Thus it follows that \eqref{eq:ThmLemma2.1Condition} holds 
	for $\Pi^{\ss}_{\aa}$-a.e.\ $\xx \in \Delta_K$ as required.
	
	Thus far we assumed $\aa \in (0,\infty)^K$. To account for the general case $\aa \in [0,\infty)^K$ we note that both sides of \eqref{eq:NonNeutralTransitionDensityWF} are weakly continuous in $\aa$, exactly as argued in \citet[Thms~3.1 and 6.2]{BarbourEthierGriffiths}.
\end{proof}

Gamma duality yields a Dirichlet-mixture representation of the transition function, where the tilted Dirichlet kernels $\Pi^{\ss}_{\aa+\lv}$ and the intractable weights $b^{\ss}_{\lv}(t;\xx)$ in \eqref{eq:BEG_transition_expansion} are replaced by standard Dirichlet kernels $\Pi^{\oo}_{\aa+\lv}$ and by the weights $b_{\lv}(t;\xx)$ induced by the dual with rates \eqref{eq:TransitionRatesDual}, respectively.
Consequently, simulation from the transition law reduces to computing the weights $b_{\lv}(t;\xx)$ and sampling from Dirichlet distributions.

The results in Theorems~\ref{thm:TransitionDensityWF}   and in Proposition \ref{thm:StationaryDistDual} yield, as a by-product practical and efficient procedures for sampling from (i) the non-neutral WF transition law and (ii) the stationary laws of the dual and of the WF diffusion. We detail these in the following two sections and then briefly discuss their performance in Section \ref{sec:stationary_numerics}.

\subsection{Sampling from the transition law}\label{subsec:algorithmsWFDiff}

Equation~\eqref{eq:NonNeutralTransitionDensityWF} leads to a straightforward mixture-based sampling strategy. Given a starting state $\xx\in\Delta_K$ and time $t>0$:
\begin{itemize}
	\item draw an index $\LL\in\ZZ$ from the entrance-law weights $\{b_{\lv}(t;\xx)\}_{\lv\in\ZZ}$;
	\item given $\LL = \lv$, draw the arrival state $\YY\sim\textrm{Dir}(\aa+\lv)$.
\end{itemize}
The main computational task is therefore to approximate the weights $b_{\lv}(t;\xx)$, which do not admit a closed-form expression easily.

For moderate and large time increments, we achieve this through Monte Carlo. Specifically, we simulate independent realisations of the dual process $(\NN_s)_{0\le s\le t}$ via a standard Gillespie scheme \citep{gillespie1977exact}, initialised from a state $\NN_0=\nn^{(0)}$ with $|\nn^{(0)}|$ large and $\nn^{(0)}/|\nn^{(0)}|\approx\xx$ so as to approximate the entrance boundary in \eqref{eq:StartFromInfinity}. Empirical frequencies of $\NN_t$ then provide estimates of $\{b_{\lv}(t;\xx)\}$, from which we sample $\lv$.

For very small $t$, accurately mimicking the entrance boundary requires $|\NN_0|$ of larger order than $t^{-1}$, and the quadratic death rates can make direct Gillespie simulation incur numerical overflows. In Section \ref{sec:stationary_numerics} we set a default threshold $\varepsilon=0.05$, and use the CLT from Theorem~\ref{thm:CLT} for $t\le\varepsilon$. In practice, we first generate a plausible total count $N\approx|\NN_t|$ using the Gaussian approximation, and then draw the composition of $\NN_t$ by multinomial splitting using a suitably perturbed parameter $\tilde{\xx}$. 

\begin{algorithm}[t]
	\DontPrintSemicolon
	\KwIn{Initial state $\xx \in \Delta_K$, time $t > 0$, time threshold $\varepsilon>0$, mutation parameters $\aa\in [0,\infty)^K$, selection parameters $\ss\in[0,\infty)^K$}
	\If{$t > \varepsilon$}{
		Approximate $\{b_{\lv}(t;\xx)\}_{\lv\in\ZZ}$ via Gillespie simulation of the dual started from $\boldsymbol{n}^{(0)}$ with $|\boldsymbol{n}^{(0)}|$ large and $\boldsymbol{n}^{(0)} / |\boldsymbol{n}^{(0)}| \approx \boldsymbol{x}$\;
		Sample $\LL \sim \{b_{\mm}(t,\xx)\}_{\mm\in\ZZ}$\;
	}
	\Else{
		Sample $Z\sim \mathcal{N}(0,1)$ and set $N = \max\Big\{0,\ \big\lfloor v(t)+\sqrt{v(t)/3}\,Z \big\rfloor \Big\}$  with $v(t)$ as in \eqref{eq:v(t)}\;
		Sample $\LL \mid N \sim \MN(\cdot;\, N,\tilde\xx)$ where $\tilde\xx=\xx(1+(t/2)(\ss-\tilde s))$\;
	}
	Given $\LL=\lv$, sample $\YY \sim \textnormal{Dir}(\aa + \lv)$\;
	\Return{$\YY$}\;
	\caption{\footnotesize Simulation from the transition function of a non-neutral WF diffusion.}\label{alg:WFDiff}
\end{algorithm}

When applying the Gaussian approximation for $\left|\NN_t\right|$, accurately reconstructing $\NN_t$ requires both information on the initial configuration $\xx$, as well as $\ss$. A multinomial sampling step $\NN_t \sim \textnormal{MN}(\left|\NN_t\right|, \xx)$ only takes into account $|\ss|$ and thus fails to properly account for the influence of $\ss$ on each individual co-ordinate. Instead, observe that
\begin{equation}\nonumber
	\left(\mathcal{A}_K f_i\right)\left(\boldsymbol{m}\right) \overset{\left|\mm\right|\rightarrow\infty}{\longrightarrow} \frac{x_i}{2}\left(\sigma_i - \tilde{s}\right), \qquad \tilde{s}:= \langle\ss,\xx\rangle
\end{equation} 
for functions $f_i(\mm) = m_i/|\mm|$. From this we deduce that
\begin{equation}\nonumber
	\mathbb{E}\left[f_i(\boldsymbol{N}_t) | \boldsymbol{N}_0 = \boldsymbol{m}\right] \approx f_i(\boldsymbol{m}) + t\mathcal{A}_Kf_i(\boldsymbol{N}_t).
\end{equation} 
Thus in Algorithm \ref{alg:WFDiff} we sample $\NN_t \sim \textnormal{MN}(\left|\NN_t\right|, \tilde\xx)$, for $\tilde\xx = \xx(\boldsymbol{1}+(t/2)(\ss-\tilde{\boldsymbol{s}}))$ where the multiplication should be interpreted as being element-wise, and we define $\tilde{\boldsymbol{s}} := \tilde{s}\boldsymbol{1}$.

The sampling scheme described in this section is summarized in Algorithm~\ref{alg:WFDiff}. Our approach applies for any fixed number of alleles $K$ and arbitrary genic selection vector $\ss$, therefore covering regimes where existing exact or numerical methods are limited, as explained above. 

\subsection{Sampling from the stationary distribution}\label{subsec:algorithmsStationaryDual}

When applied for large $t$, Algorithm~\ref{alg:WFDiff} approximates draws from the WF stationary law by sampling from the transition function. However, the dual process also yields a direct route to stationarity. From \eqref{eq:NonNeutralTransitionDensityWF}, letting $t\to\infty$, the mixture weights converge to the stationary law $\pi$ of the dual process, yielding a mixture representation of \eqref{eq:StationaryWF} as
\begin{equation}\label{eq:Stationary_mixture_generic}
	\pi_{\aa}^{\ss}(\yy)
	=
	\sum_{\lv\in\ZZ}\pi(\lv)\,\pi^{\oo}_{\aa+\lv}(\yy)
\end{equation}
when $\aa \in (0,\infty)^K$. Accordingly, sampling from the stationary distribution $\Pi_{\aa}^{\ss}$ reduces to drawing $\LL\sim\pi$ with mass function \eqref{eq:DualStationaryDistribution} and then sampling
$\YY \mid \LL=\lv \sim \mathrm{Dir}(\aa+\lv)$.

In principle, sampling from \eqref{eq:DualStationaryDistribution} involves evaluating the normalising constant $\Phi^{(K)}_{2}(\aa;|\aa|;\ss)$ in \eqref{eq:normalisingConstant}. In the Appendix we describe two exact accept--reject schemes that avoid evaluating $\Phi^{(K)}_{2}$ altogether. Despite their theoretical appeal, these schemes exhibit low acceptance rates under strong selection and/or increasing $K$. Hence, we  propose here a new, approximate but substantially more efficient sampling scheme based on truncation, whose cost is essentially insensitive to the selection coefficients and increases only mildly with $K$. We also provide an explicit total-variation bound that allows choosing the truncation level for a user-specified tolerance.


Regroup \eqref{eq:normalisingConstant} by total degree $m=|\lv|$, to get
\begin{equation}\label{eq:Phi2_degree_grouping_alg}
	\Phi^{(K)}_{2}\!\left(\aa;|\aa|;\ss\right)
	=
	\sum_{m=0}^{\infty}\frac{a_m}{|\aa|_{(m)}},
	\quad
	a_m
	:=
	\sum_{\lv\in\ZZ:\, |\lv|=m}
	\frac{\aa_{(\lv)}}{\lv!}\,\ss^{\lv},\quad m\ge1,
\end{equation}
with $a_0:=1$.
Thus the (unnormalised) total-degree weights are $T_m:=a_m/|\aa|_{(m)}$ and the
scheme samples $M:=|\LL|$ from the pmf defined by $\{T_m\}_{m\ge0}$, then draws $\LL$ conditional on $M$,
and finally draws $\YY\mid \LL$ from the corresponding Dirichlet component.

The coefficients $a_m$ admit an equivalent representation in terms of partitions of $m$, made explicit by the following lemma.

\begin{lemma}\label{lemma:Phi2}
	For $K\in\mathbb{N}$, $\aa\in(0,\infty)^K$, and $\ss\in[0,\infty)^K$,
	\begin{align}\label{eq:Phi2Reduction_alg}
		\Phi^{(K)}_{2}\!\left(\aa;|\aa|;\ss\right)
		=
		1+\sum_{m=1}^{\infty}\frac{1}{|\aa|_{(m)}}\sum_{\kk\in\mathbb{K}_m}
		\prod_{j=1}^{m}\frac{1}{k_j!}
		\left(\frac{\sum_{i=1}^{K}\alpha_i\sigma_i^{\,j}}{j}\right)^{k_j},
	\end{align}
	where $\mathbb{K}_m
	:=
	\{\kk=(k_1,\dots,k_m)\in\mathbb{Z}_+^{\,m}:\ \sum_{j=1}^m j k_j = m\}$.
\end{lemma}

\begin{proof}
	By Theorem 1 in \cite{VanLaarhovenKalker88} we have the following identity for the Lauricella function $F_{D}\left(a; \bb; c; \xx\right)$
\begin{equation}\nonumber
\begin{aligned}
		F_{D}\left(a; \bb; c; \xx\right) = 1 + \sum_{m=1}^{\infty}\frac{a_{\left(m\right)}}{c_{\left(m\right)}}\sum_{\kk\in\mathbb{K}_m}\prod_{j=1}^{m}\frac{1}{k_j!}\left(\frac{\sum_{i=1}^{K}b_ix_i^j}{j}\right)^{k_j}
\end{aligned}
\end{equation} 
	where $\mathbb{K}_m := \left\{\kk \in \ZZ : \sum_{j=1}^{K} jk_j = m\right\}$. 
Equation \eqref{eq:Phi2Reduction_alg} follows by applying the confluence limit to the identity for $F_D$. Specifically, using the relation
\begin{equation}\nonumber
\Phi^{(K)}_{2}\left(\aa; |\aa|; \ss\right) = \lim_{a\to\infty} F_{D}\left(a; \aa; |\aa|; \ss/a\right),
\end{equation}
and noting that $\lim_{a\to\infty} a_{(m)}/a^m = 1$, we substitute $x_i = \sigma_i/a$ into the expansion and take the limit term-wise, 
which yields \eqref{eq:Phi2Reduction_alg} as required.
\end{proof}

Comparing \eqref{eq:Phi2_degree_grouping_alg} and \eqref{eq:Phi2Reduction_alg} yields, for $m\ge1$,
\begin{equation}\nonumber
	a_m
	=
	\sum_{\kk\in\mathbb{K}_m}
	\prod_{j=1}^{m}\frac{1}{k_j!}
	\left(\frac{\sum_{i=1}^{K}\alpha_i\sigma_i^{\,j}}{j}\right)^{k_j},
\end{equation} 
(cf.\ \citealp[eq.\ (14)]{VanLaarhovenKalker88}).
For computation, it is convenient to exploit the following recurrence.

\begin{lemma}\label{lem:am_recursion}
	Let $a_m$ be as in \eqref{eq:Phi2_degree_grouping_alg} and define
	$S_j:=\sum_{i=1}^K \alpha_i\sigma_i^{\,j}$ for $j\ge1$. Then
	\begin{equation}\label{eq:a_m_recurrence_alg}
		a_0=1,
		\qquad
		a_m
		=
		\frac{1}{m}\sum_{j=1}^{m} S_j\,a_{m-j},
		\qquad m\ge 1.
	\end{equation}
\end{lemma}
\begin{proof}
    Define the generating function $A(z) := \sum_{m\ge 0} a_m z^m$. By the definition of $a_m$ and the binomial series identity, we have the product representation $A(z) = \prod_{i=1}^K (1-\sigma_i z)^{-\alpha_i}$. Taking logarithms and expanding via $-\log(1-x) = \sum_{j\ge 1} x^j/j$ yields
    \[
    \log A(z) = \sum_{i=1}^K \alpha_i \sum_{j\ge 1} \frac{(\sigma_i z)^j}{j} = \sum_{j\ge 1} \frac{z^j}{j} \left( \sum_{i=1}^K \alpha_i \sigma_i^j \right) = \sum_{j\ge 1} \frac{S_j}{j} z^j.
    \] 
    Differentiating both sides with respect to $z$ gives $A'(z) = A(z) \sum_{j\ge 1} S_j z^{j-1}$.
     Equating the coefficients of $z^{m-1}$ on both sides of this identity results in $    m a_m = \sum_{j=1}^m S_j a_{m-j}$ which gives the result.
\end{proof}

Now let $M$ have mass function $\PP(M=m)\propto T_m$ on $\mathbb{Z}_+$ and define the conditional law
\begin{equation}\nonumber
	\PP(\LL=\lv\mid M=m)\ \propto\ \frac{\aa_{(\lv)}}{\lv!}\,\ss^{\lv},
	\qquad |\lv|=m.
\end{equation} 
In practice we truncate at $M_{\max}$ and renormalise on $\{0,\dots,M_{\max}\}$.
Denote by $\Pi_{\aa}^{\ss,M_{\max}}$ the induced law of $\YY$ on $\Delta_K$.
The next proposition gives an explicit total-variation bound.

\begin{proposition}\label{prop:TV_truncation_main}
	Let $\Pi_{\aa}^{\ss}$ denote the target stationary WF distribution on $\Delta_K$
	with density $\pi_{\aa}^{\ss}$ in \eqref{eq:StationaryWF}, and let $\Pi_{\aa}^{\ss,M_{\max}}$
	be the law induced by the truncated total-degree scheme. Then
	\begin{equation}\label{eq:TV_truncation_main}
		\bigl\|\Pi_{\aa}^{\ss}-\Pi_{\aa}^{\ss,M_{\max}}\bigr\|_{\mathrm{TV}}
		\;\le\;
		2e^{|\ss|}\Big(\frac{e\,|\ss|}{M_{\max}+1}\Big)^{M_{\max}+1},
		\qquad
		|\ss|:=\sum_{i=1}^K \sigma_i,
	\end{equation}
	where $\|\cdot\|_{\mathrm{TV}}$ denotes total variation norm.
\end{proposition}

\begin{proof}
	Let $\Pi_{\aa}^{\ss}$ be the law induced by the infinite total-degree scheme, and
	$\Pi_{\aa}^{\ss,M_{\max}}$ the law induced by truncating $M$ to $\{0,\dots,M_{\max}\}$
	and renormalising. Writing both as mixtures over the common conditional laws
	$\Pi(\cdot\mid M=m)$ gives
	\[
	\Pi_{\aa}^{\ss}(\cdot)=\sum_{m\ge0}p_m\,\Pi(\cdot\mid M=m),
	\qquad
	\Pi_{\aa}^{\ss,M_{\max}}(\cdot)=\sum_{m=0}^{M_{\max}}\tilde p_m\,\Pi(\cdot\mid M=m),
	\]
	with $p_m=T_m/\sum_{r\ge0}T_r$ and $\tilde p_m=p_m/\sum_{r=0}^{M_{\max}}p_r$ for
	$m\le M_{\max}$. Hence, for any measurable $A\subseteq\Delta_K$,
	\[
	\big|\Pi_{\aa}^{\ss}(A)-\Pi_{\aa}^{\ss,M_{\max}}(A)\big|
	\le
	\sum_{m=0}^{M_{\max}}|p_m-\tilde p_m|+\sum_{m>M_{\max}}p_m
	=
	2\sum_{m>M_{\max}}p_m,
	\]
	so
	\begin{equation}\label{eq:TVboundPartial}
		\bigl\|\Pi_{\aa}^{\ss}-\Pi_{\aa}^{\ss,M_{\max}}\bigr\|_{\mathrm{TV}}
		\le
		2\sum_{m>M_{\max}}p_m
		=
		2\,\frac{\sum_{m>M_{\max}}T_m}{\sum_{r\ge0}T_r}
		\le
		2\sum_{m>M_{\max}}T_m,
	\end{equation}
	using $\sum_{r\ge0}T_r=\Phi^{(K)}_{2}(\aa;|\aa|;\ss)\ge 1$.
	
	To bound the discarded mass, fix $m\ge1$ and note that for $|\lv|=m$,
	$\aa_{(\lv)}\le(|\aa|)_{(m)}$, hence
	\[
	a_m=\sum_{|\lv|=m}\frac{\aa_{(\lv)}}{\lv!}\,\ss^{\lv}
	\le (|\aa|)_{(m)}\sum_{|\lv|=m}\frac{\ss^{\lv}}{\lv!}
	=
	(|\aa|)_{(m)}\,\frac{\left|\ss\right|^{m}}{m!},
	\]
	where the last equality is the multinomial identity. Therefore
	$T_m=a_m/|\aa|_{(m)}\le \left|\ss\right|^{m}/m!$, and
	\[
	\sum_{m>M_{\max}}T_m\le \sum_{m>M_{\max}}\frac{\left|\ss\right|^{m}}{m!}
	= e^{\left|\ss\right|}\PP(X>M_{\max}),\qquad X\sim\mathrm{Poisson}(\left|\ss\right|).
	\]
	Applying a Chernoff bound yields
	$\PP(X>M_{\max})\le \big(\tfrac{e\,\left|\ss\right|}{M_{\max}+1}\big)^{M_{\max}+1}$,
	which combined with \eqref{eq:TVboundPartial} gives \eqref{eq:TV_truncation_main}.
\end{proof}

\begin{algorithm}[t]
	\DontPrintSemicolon
	\KwIn{Mutation parameters $\aa$, selection parameters $\ss$, tolerance $\epsilon>0$}
	Choose $M_{\max}$ such that $2e^{|\ss|}\Big(\tfrac{e\,|\ss|}{M_{\max}+1}\Big)^{M_{\max}+1}\le \epsilon$\;
	Compute $\{a_m\}_{m=0}^{M_{\max}}$ via \eqref{eq:a_m_recurrence_alg} and set $T_m:=a_m/|\aa|_{(m)}$\;
	Sample $M\in\{0,\dots,M_{\max}\}$ with $\PP(M=m)\propto T_m$\;
	Given $M=m$, sample $\LL$ on $\{\lv\in\ZZ:\,|\LL|=m\}$ with $\PP(\LL=\lv\mid M=m)\propto \aa_{(\lv)}\ss^{\lv}/\lv!$\;
	Given $\LL=\lv$, sample $\YY \sim \textnormal{Dir}(\aa+\lv)$\;
	\Return{$\YY$}\;
	\caption{\footnotesize Approximate simulation from the stationary distribution \eqref{eq:StationaryWF} with explicit error control.}\label{alg:StationaryWF}
\end{algorithm}

The above simulation scheme with explicit error control is summarized in Algorithm~\ref{alg:StationaryWF}. 

\subsection{Numerical illustration}\label{sec:stationary_numerics}

\begin{figure}[t]
	\centering\includegraphics[width=.42\linewidth]{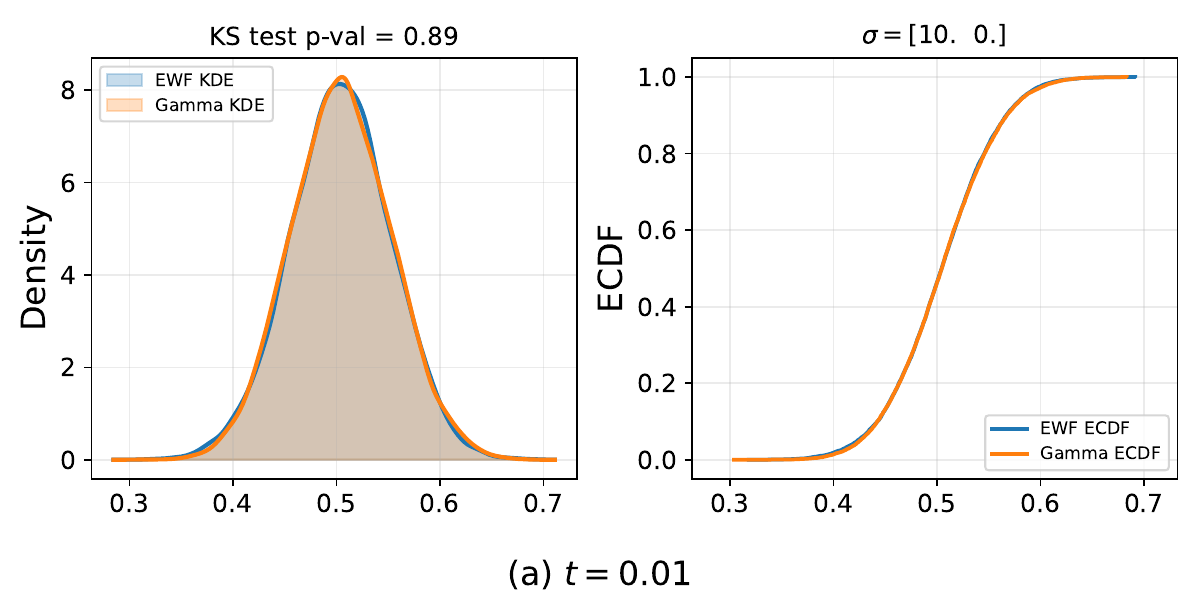}
	\includegraphics[width=.42\linewidth]{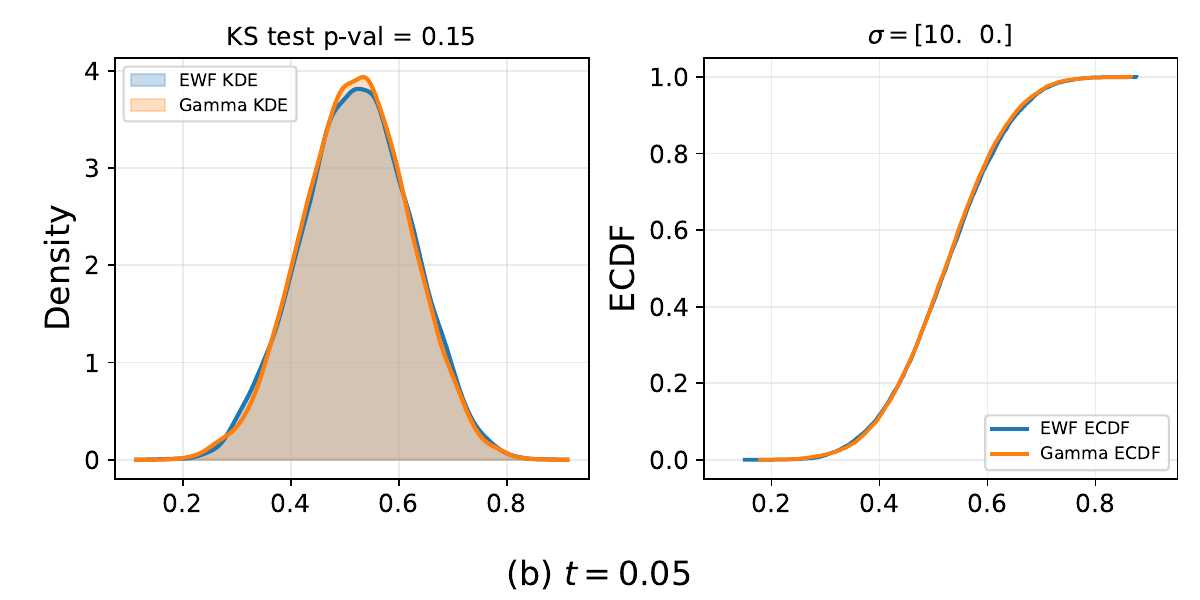}\\
	\centering\includegraphics[width=.42\linewidth]{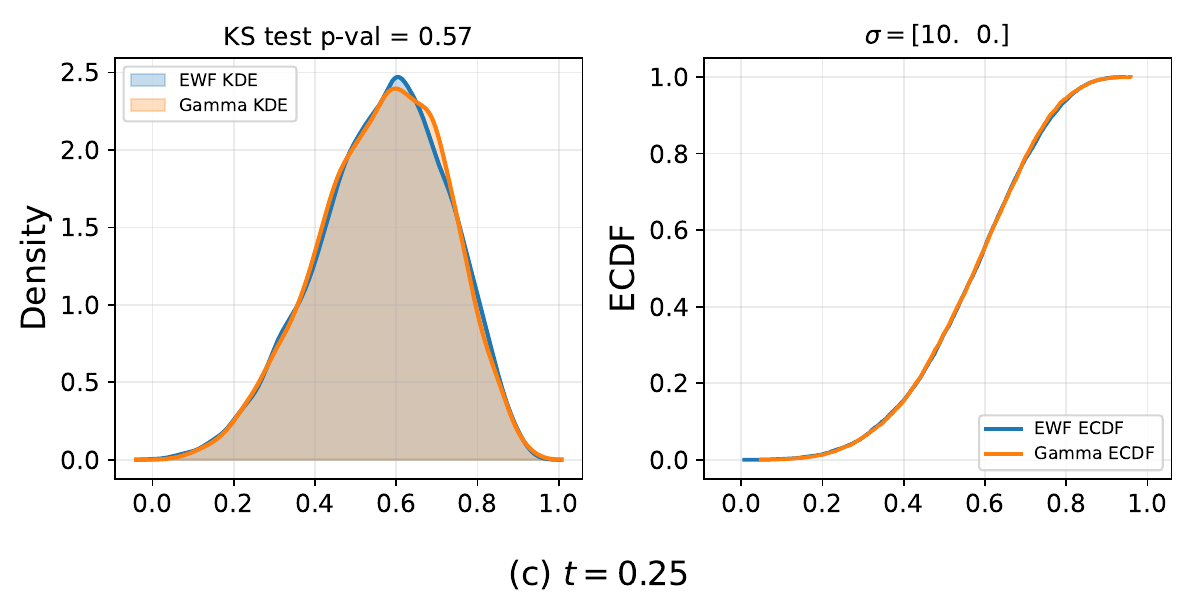}
	\includegraphics[width=.42\linewidth]{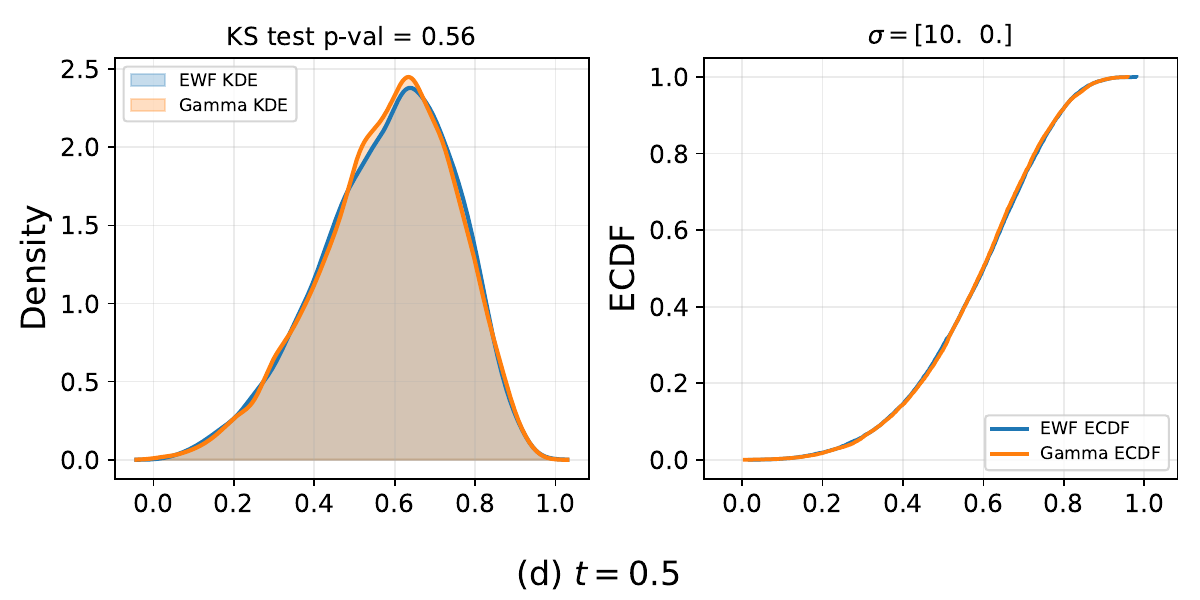}
	\caption{\footnotesize Comparison of the histograms (with corresponding kernel smoothing) and empirical CDFs between the samples generated using the dual-based sampler in Alg.\ \ref{alg:WFDiff} and \texttt{EWF} (blue) for (a) $t = 0.01$, (b) $t = 0.05$, (c) $t = 0.25$, and (d) $t = 0.5$. Across all simulations we set $\aa = (2,5)$, $\xx=(0.5, 0.5)$, and $\varepsilon = 0.05$, thus the top two panels make use of the Gaussian approximation from Theorem \ref{thm:CLT}, whilst the bottom two panels apply a Gillespie algorithm. We further run a two sample Kolmogorov--Smirnov test of goodness of fit to check for statistically significant differences between the two samples. All reported $p$-values are above the 0.05 level of significance across all simulation scenarios.}
	\label{fig:GammaDualityvEWF}
\end{figure}

\begin{figure}[t]
	\centering\includegraphics[width=.8\linewidth]{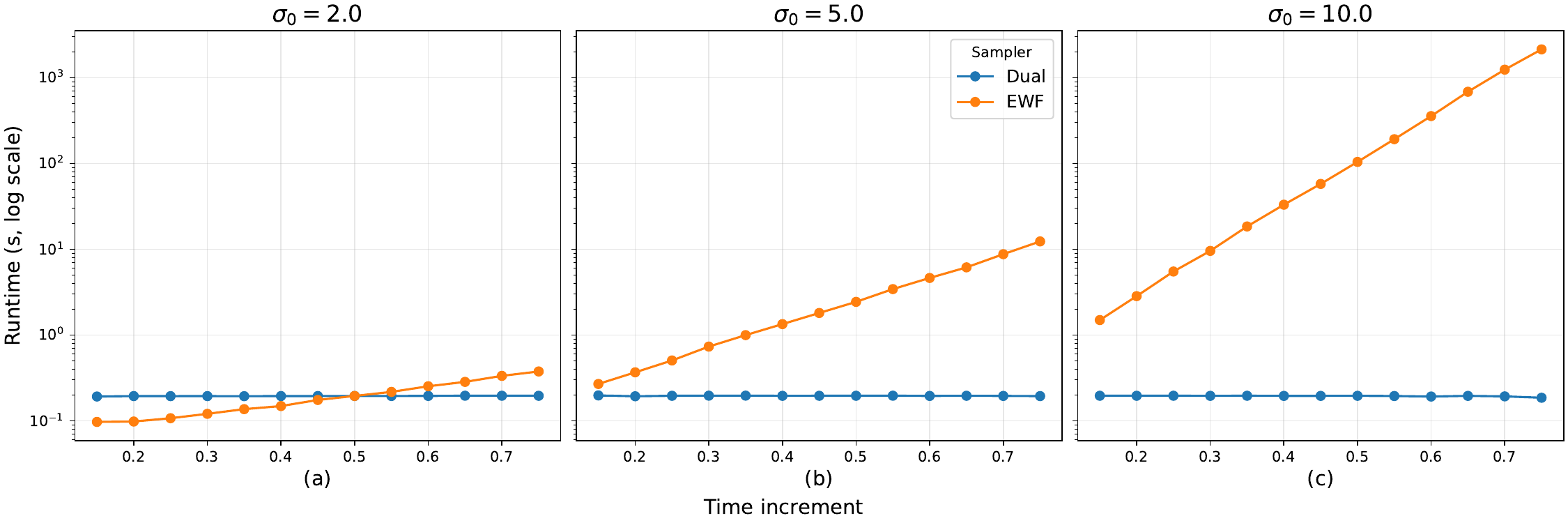}
	\caption{\footnotesize Log-scale runtime comparison in the one-dimensional ($K=2$) WF case as the selection coefficient increases, between the proposed dual-based sampler and the \texttt{EWF} sampler \citep{SantEWF}. Across all scenarios, we set $\aa = (2,5)$, $\xx = (0.5, 0.5)$, and consider equally space time increments of 0.05 from 0.1 to 0.75. Dual-based runtimes are essentially constant across selection values and time increments.}
	\label{fig:GammaDualityvEWFRuntime}
\end{figure}

We first illustrate the application of Algorithm \ref{alg:WFDiff} for sampling from \eqref{eq:NonNeutralTransitionDensityWF}.
Fig.~\ref{fig:GammaDualityvEWF} shows that our method closely matches the exact sampler \texttt{EWF} \citep{SantEWF} in dimension $K=2$, in both regimes when the Gaussian limit in Theorem \ref{thm:CLT} is used (top panels) and those in which the Gillespie algorithm is used (bottom panels). Fig.\ \ref{fig:GammaDualityvEWFRuntime} compares runtimes for generating 10,000 samples using both samplers as the selection coefficient increases. Notably, the dual-based runtimes remain essentially constant across selection values and time increments. To our knowledge, no alternative method is available through standard software for analogous runtime comparisons when $K\ge 3$. 

\begin{figure}[t]
	\centering
	\includegraphics[width=.8\linewidth]{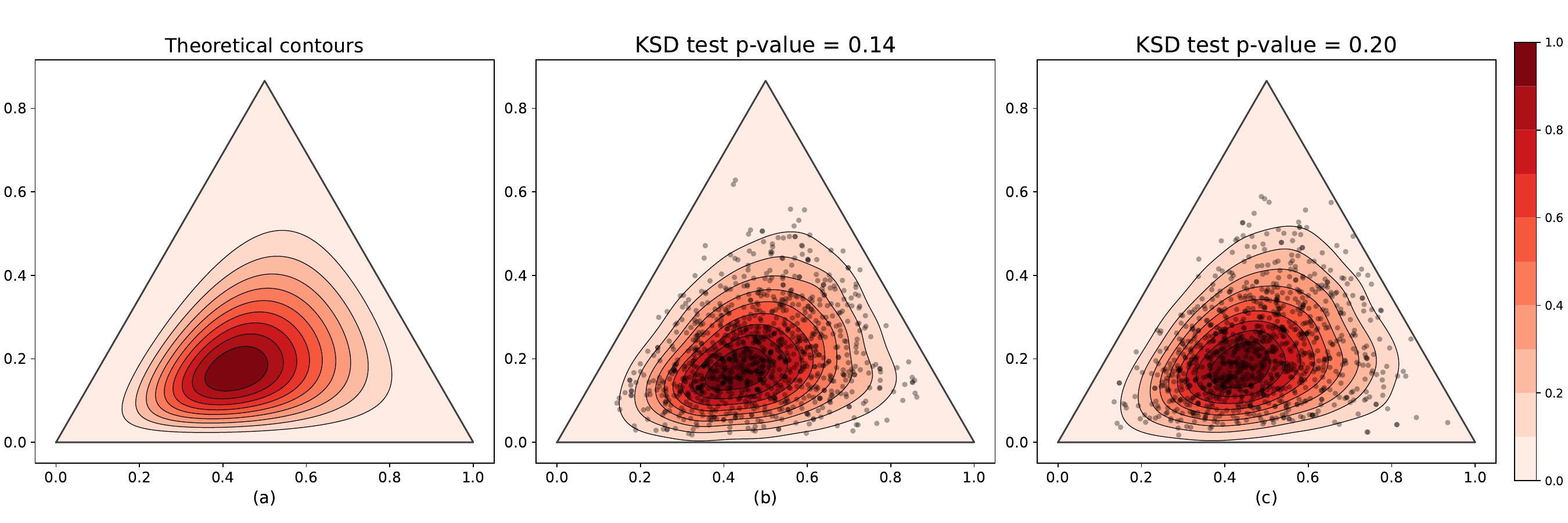}
	\caption{\footnotesize Stationary distribution \eqref{eq:StationaryWF} for $K=3$ with $\boldsymbol{\sigma}=(10,0,2)$ and $\boldsymbol{\alpha}=(2,5,3)$, shown via the barycentric embedding of the $2$-simplex in $\mathbb{R}^2$. Panel (a) shows the theoretical density; panel (b) shows the empirical contour plot from Algorithm~\ref{alg:WFDiff} at $t=50$; panel (c) shows the empirical contour plot from Algorithm~\ref{alg:StationaryWF} based on $10{,}000$ samples. In each case, the KSD $p$-value shows no statistical evidence to distinguish the empirical output from the target stationary distribution.}
	\label{fig:GammaStationaryK3}
\end{figure}

Fig.\ \ref{fig:GammaStationaryK3} shows, for $K=3$ with $\boldsymbol{\sigma}=(10,0,2)$ and $\boldsymbol{\alpha}=(2,5,3)$, the theoretical stationary density \eqref{eq:StationaryWF} in panel (a) and an empirical contour plot from the dual-based sampler (Algorithm~\ref{alg:WFDiff}) at $t=50$ in panel (b), both rendered via the barycentric embedding of the $2$-simplex in $\mathbb{R}^2$. A Kernel Stein discrepancy (KSD) test \citep{chwialkowski2016kernel} finds no evidence to distinguish the empirical output from the target stationary distribution, as indicated by the large $p$-value. Table~\ref{tab:ksd_results} (column 2) repeats this assessment as the number of alleles increases, again indicating statistical agreement with the theoretical distribution.

We now discuss the performance of Algorithm \ref{alg:StationaryWF}. Fig.\ \ref{fig:GammaStationaryK3} (panel (c)) shows, for $K=3$, an empirical contour plot obtained from $10{,}000$ samples, with the KSD test once again finding no evidence to reject equality in distribution to the target.
Table~\ref{tab:ksd_results} (last columns) indicates that this agreement persists as $K$ increases.

\begin{figure}[t]
	\centering\includegraphics[width=.8\linewidth]{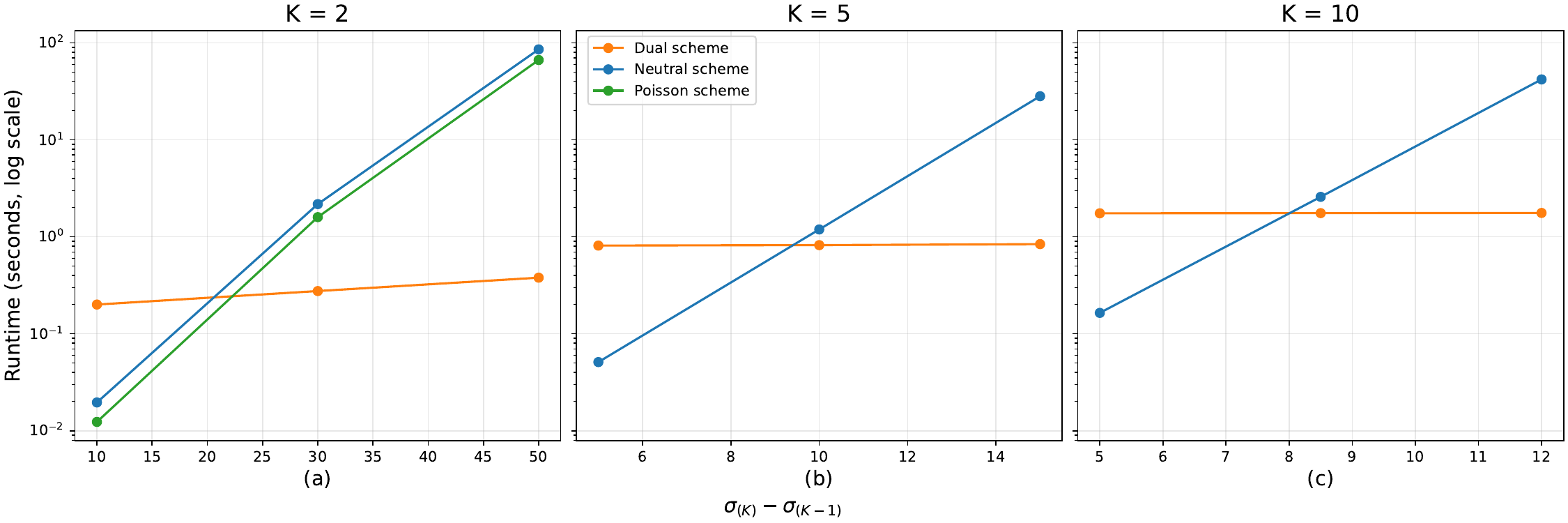}
	\caption{\footnotesize Log-scale runtime comparison for sampling from the stationary law as the selection gap
		$\sigma_{(K)}-\sigma_{(K-1)}$ increases, shown for $K=2,5,10$.
		Displayed are the dual-based total-degree scheme (Algorithm~\ref{alg:StationaryWF}), the neutral accept--reject scheme, and the product-Poisson accept--reject scheme for the dual weights (shown only for $K=2$, as for $K\ge 3$ its runtimes exceed the plotted range).}
	\label{fig:RejectionSampler}
\end{figure}

Fig.\ \ref{fig:RejectionSampler} compares runtimes for generating 10,000 samples when $K=2,5,10$ as the selection gap
$\sigma_{(K)}-\sigma_{(K-1)}$ between the two largest selection coefficients increases. For $K=2$ this gap fully determines the exponential tilt, while for $K\ge 3$ it provides a convenient one-dimensional proxy for selection strength (with the remaining gaps held fixed), under which accept--reject efficiency typically deteriorates.
The figure reports the dual-based total-degree scheme (Algorithm~\ref{alg:StationaryWF}) alongside the two exact accept--reject schemes detailed in the Appendix: the sampler based on neutral Dirichlet proposals, and the one using product-Poisson proposals.
Across the plotted regimes, the total-degree scheme exhibits essentially constant runtime as the selection gap increases and very mild growth with $K$.
By contrast, the runtimes for the exact schemes increase rapidly as this gap increases, with the product-Poisson sampler becoming prohibitively slow already for $K=5$.

Additional numerical experiments in the Appendix investigate: (i) accuracy checks for the evaluation of $\Phi^{(K)}_{2}(\aa;|\aa|;\ss)$; (ii) agreement between the Gaussian approximation in Theorem~\ref{thm:CLT} and Algorithm \ref{alg:WFDiff} with a suitably large initialisation; and (iii) sensitivity of the Monte Carlo estimates of the weights $\{b_{\lv}(t;\xx)\}_{\lv\in\ZZ}$ (via the Gillespie simulation) to the initialisation of the birth-and-death process.

\begin{table}[t]
	\centering
	\begin{tabular}{r|cc|cc}
		& \multicolumn{2}{c|}{$t = 50$} & \multicolumn{2}{c}{Stationary} \\
		$K$ & $\mathrm{KSD}$ statistic & $p$-value & $\mathrm{KSD}$ statistic & $p$-value \\
		\hline
		4  & $-7.8915\times 10^{-3}$  & $0.63$ & $7.1252\times 10^{-3}$  & $0.23$\\
		5  & $2.4274\times 10^{-2}$ & $0.22$ & $-3.6816\times 10^{-2}$  & $0.97$ \\
		10 & $-2.9051\times 10^{-1}$ & $1.00$ & $7.8087\times 10^{-3}$  & $0.27$\\
		20 & $-1.0514\times 10^{-3}$ & $1.00$ & $-1.6424\times 10^{-4}$  & $0.58$\\
		\hline
	\end{tabular}
	\caption{\footnotesize Kernel Stein Discrepancy (KSD) test statistics and corresponding $p$-value for increasing dimension $K$ when using both Algorithm \ref{alg:WFDiff} (and setting $t = 50$), as well as Algorithm \ref{alg:StationaryWF}.}
	\label{tab:ksd_results}
\end{table}

\section{Discussion}\label{Discussion}

We introduced a gamma duality framework for a multi-allelic WF diffusion with parent-independent mutation and genic selection. The construction yields a multivariate birth-and-death dual with explicit infinitesimal rates, an identified stationary distribution, and a small-time description, and it leads to a Dirichlet-mixture representation of the transition function with mixing weights governed by the dual started from an entrance boundary. This representation provides an efficient route to computation, including regimes where selection is strong and the number of alleles is arbitrary.

Several directions remain open. First, our analysis is restricted to genic selection, and the particular exponential augmentation of polynomial duality used here does not seem to directly extend to frequency-dependent selection. It would be mathematically and biologically valuable to determine whether related duality constructions exist in that setting, and whether they can again be made computationally tractable. Second, while the present algorithms evolve the dual via a Gillespie-type simulation, it is natural to ask whether a spectral characterisation can be recovered for the new dual, which could in turn lead to fully exact transition samplers without Monte Carlo approximation of the weights. Third, extending the framework to the WF bridge diffusion would broaden its usefulness for data-augmentation and MCMC-based inference, but doing so in the non-neutral setting requires additional ideas to control conditioning and maintain tractability.	Finally, the new dual process has the potential to offer new insights into the genealogy of a sample under selection. Currently we do not have an intuitive interpretation for the jump rates in \eqref{eq:TransitionRatesDual}. Contrary to the ancestral selection graph (see, e.g., \cite{baa:wak:2018} for review), it is the \emph{fitter}, not weaker, types that experience greater birth rates, so the standard interpretation of a birth event as tracking `potential parents' does not apply.





\section*{Acknowledgements}

	PJ was supported in part by the Engineering and Physical Sciences Research Council (EPSRC) via ProbAI: A Hub for the Mathematical and Computational Foundations of Probabilistic AI (EP/Y028783/1).
	JS was supported by the Italian Ministry of University and Research (MUR) under the PRIN 2022 programme, Grant no. 2022CLTYP4.
	MR was supported by the European Union -- NextGenerationEU, PRIN 2022 PNRR, Grant no. P2022H5WZ9.
	
	


\bibliographystyle{apalike}
\bibliography{bibfile}

\appendix
\section*{Appendix}

\subsection*{Exact sampling from the stationary distribution}

We describe two exact accept--reject schemes to sample from $\Pi_{\aa}^{\ss}$ that bypass evaluation of the normalising constant $\Phi^{(K)}_{2}(\aa;|\aa|;\ss)$.

The WF stationary density has exponential-tilt form $\pi_{\aa}^{\ss}(\yy)\propto \pi^{\oo}_{\aa}(\yy)\exp\{\langle \ss,\yy\rangle\}$ on $\Delta_K$. Since $\langle \ss,\yy\rangle\le \sigma_{(K)}:=\max_{1\le i\le K}\sigma_i$ for all $\yy\in\Delta_K$, a standard accept--reject sampler proposes $\YY\sim\mathrm{Dir}(\aa)$ (neutral model) and accepts with probability $\exp\{\langle \ss,\YY\rangle-\sigma_{(K)}\}$. An analogous strategy appears in Algorithm~1 of \cite{donnelly2001likelihoods} for simulating samples under selection. We shall henceforth refer to this sampling strategy as the neutral accept--reject scheme.

Alternatively, we sample the dual index $\LL\in\ZZ$ from its stationary mass function \eqref{eq:DualStationaryDistribution} and then draw the Dirichlet component $\YY\mid \LL=\lv\sim \mathrm{Dir}(\aa+\lv)$. To allow the first of these steps, let $q(\lv)$ be the product-Poisson($\ss$) law on $\ZZ$ with independent coordinates.  From \eqref{eq:DualStationaryDistribution},
\[
\pi(\lv)
\;\propto\;
\frac{\aa_{(\lv)}}{|\aa|_{(|\lv|)}}\,\frac{\ss^{\lv}}{\lv!}
\;\propto\;
q(\lv)\,\frac{\aa_{(\lv)}}{|\aa|_{(|\lv|)}},
\]
so an accept--reject sampler with proposal $q(\cdot)$ uses acceptance probability
\begin{equation}\nonumber
	\beta(\lv)
	=
	\frac{\aa_{(\lv)}/|\aa|_{(|\lv|)}}{\sup_{\kk\in\ZZ}\ \aa_{(\kk)}/|\aa|_{(|\kk|)}}.
\end{equation}
Using $|\lv|=\sum_{i=1}^K \ell_i$ and $(\alpha_i)_{(\ell_i)}\le (|\aa|)_{(\ell_i)}$, we have $0<\aa_{(\lv)}/|\aa|_{(|\lv|)}\le 1$ for all $\lv\in\ZZ$, with equality at $\lv=\boldsymbol{0}$. Hence the supremum in $\beta(\lv)$ equals $1$ and $\beta(\lv)=\aa_{(\lv)}/|\aa|_{(|\lv|)}$. This yields exact draws from \eqref{eq:DualStationaryDistribution} without evaluating $\Phi^{(K)}_{2}$, and the subsequent draw $\YY\mid \LL=\lv\sim\mathrm{Dir}(\aa+\lv)$ produces an exact stationary sample. We will refer to this algorithm as the product-Poisson accept--reject scheme. Nonetheless, both accept--reject schemes have acceptance probabilities that deteriorate quickly as $K$ grows and/or under strong selection, motivating the much more efficient scheme introduced in Section \ref{subsec:algorithmsStationaryDual}.

\subsection*{Computing the normalising constant $\Phi^{(K)}_{2}\left(\aa; \left|\aa\right|; \ss\right)$}\label{Performance:Phi}
No standard numerical library provides pointwise evaluation of $\Phi^{(K)}_{2}\left(\aa; \left|\aa\right|; \ss\right)$ in dimensions $K > 2$. We therefore exploit~\eqref{eq:Phi2Reduction_alg} and evaluate the resulting series term-wise, truncating once the contributions fall below some prescribed tolerance $\varepsilon$, motivated by Proposition~\ref{prop:TV_truncation_main}. 
		
\begin{figure}[t]
	\centering\includegraphics[width=0.8\linewidth]{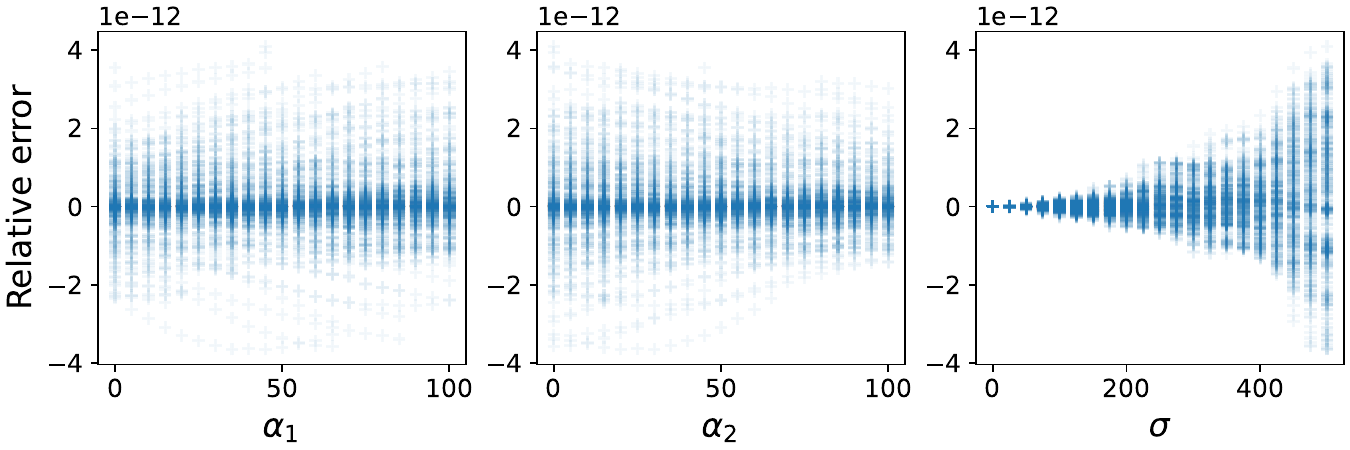}
	\caption{\footnotesize Relative error ($y$-axis) against $\alpha_1$ (left), $\alpha_2$ (middle), and $\sigma$ (right) between our implementation and that obtained when using \texttt{scipy.special.hyp1f1}. We observe that we expect a deterioration in the accuracy as $\sigma$ increases, but point out the $y$-axis scaling.}
	\label{fig:LauricellaPointwiseErrors}
\end{figure}

\begin{figure}[t]
	\centering\includegraphics[width=.5\linewidth]{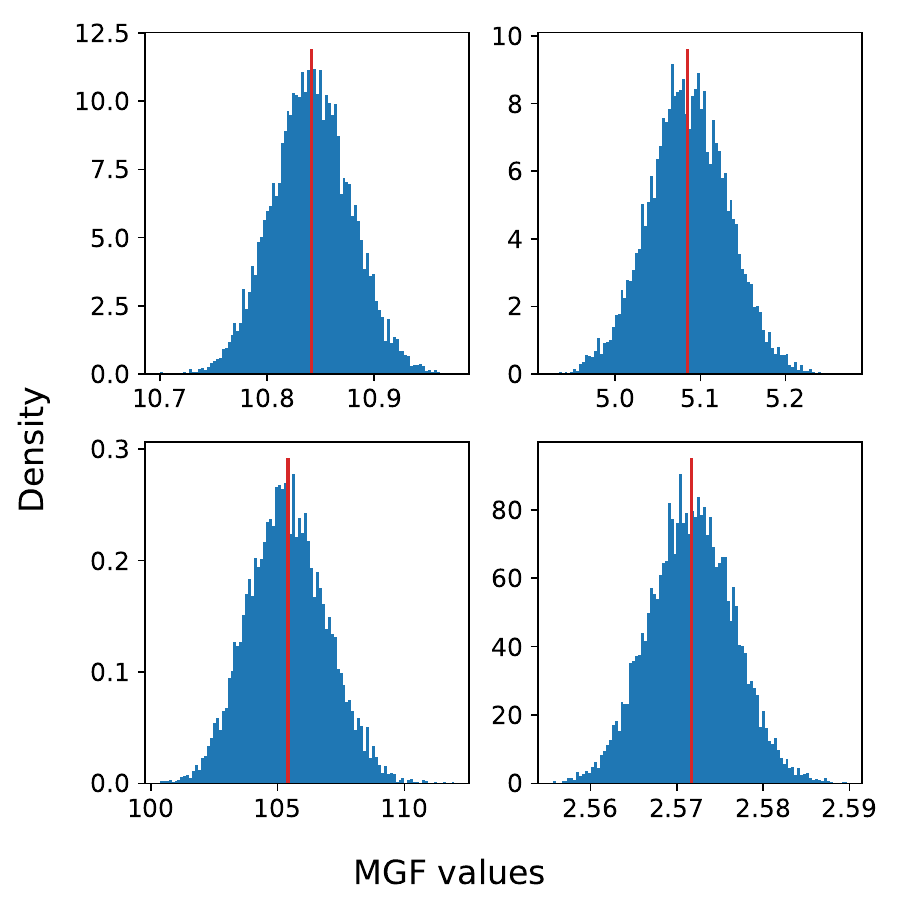}
	\caption{\footnotesize Histograms of 10,000 Monte Carlo estimates of the moment generating function for a $\textnormal{Dir}(\alpha)$ random variable with overlaid value of $\Phi^{(K)}_{2}\left(\aa; \left|\aa\right|; \ss\right)$ as computed using our implementation.}
	\label{fig:LauricellaMGFChecks}
\end{figure}

We assess the accuracy of our implementation in two ways. First, in the univariate case, we compare against standard evaluations of the confluent hypergeometric ${}_1F_1$ function (as implemented in \texttt{scipy.special.hyp1f1}). Second, in higher dimensions, we compare against Monte Carlo estimates of the moment generating function of an appropriate Dirichlet random variable, exploiting the connection in Remark~\ref{remark:DirichletConnection}. As illustrated in Figure~\ref{fig:LauricellaPointwiseErrors}, the proposed method yields accurate evaluations of $\Phi^{(K)}_{2}\left(\aa; \left|\aa\right|; \ss\right)$ across all tested parameter regimes in the univariate setting. Figure~\ref{fig:LauricellaMGFChecks} demonstrates strong agreement between our implementation and the Monte Carlo estimates in higher dimensions. 
		
\subsection*{Gaussian approximation}\label{Performance:CLT}
\begin{figure}[t]
	\centering\includegraphics[width=.7\linewidth]{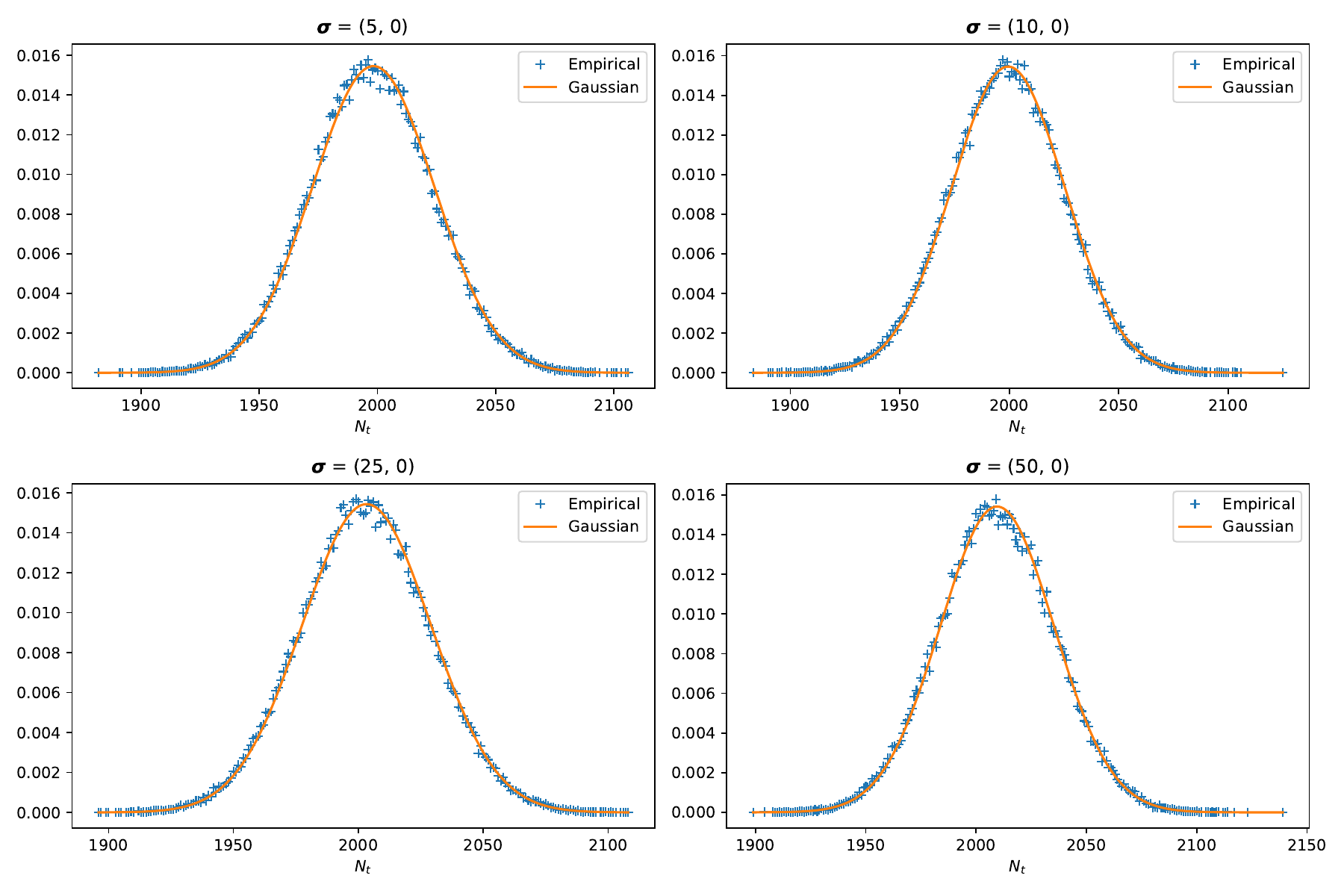}
	\caption{\footnotesize Monte Carlo estimate of the probability mass function (blue crosses) of $N_t$ for $t=0.01$ with the theoretical distribution superimposed (orange). The selection parameter $\ss$ varies as follows: 5 (top left), 10 (top right), 25 (bottom left), and 50 (bottom right).}
	\label{fig:CLT}
\end{figure}

We illustrate the convergence in Theorem~\ref{thm:CLT} by comparing the empirical distribution of $\lvert \boldsymbol{N}_t \rvert$ to the Gaussian limit obtained therein. Specifically, for $t=0.01$ we generate $10^4$ Monte Carlo samples of $\boldsymbol{N}_t$, compute the corresponding empirical probability mass function of $\lvert \boldsymbol{N}_t \rvert$, and superimpose the corresponding Gaussian distribution in Figure~\ref{fig:CLT}. We consider increasing values of the selection parameter, which induces the expected rightward shift in the mean, and observe close agreement across all regimes. In light of the sensitivity analysis in Figure~\ref{fig:InitPop}, we initialise the Gillespie algorithm with total population size equal to $10^7$, to ensure reliable estimation in this small-time regime. 
		
\subsection*{Initial population size for the Gillespie algorithm}
\begin{figure}[t]
	\centering\includegraphics[width=.9\linewidth]{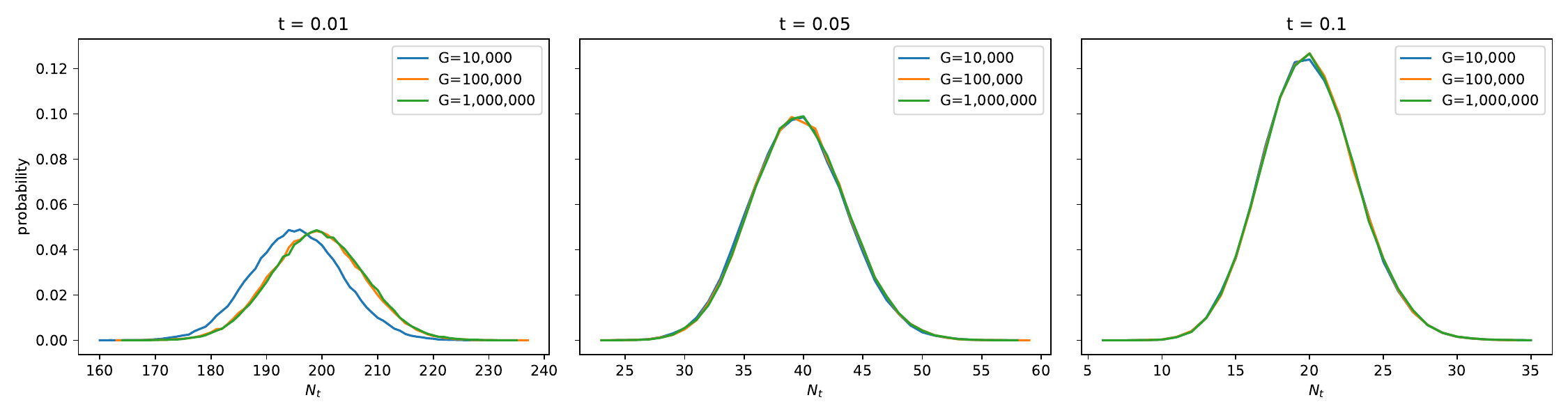}
	\caption{\footnotesize Monte Carlo estimate of the probability mass function (blue crosses) of $N_t$ for $t=0.01$ with the theoretical distribution superimposed (orange). The selection parameter $\ss$ varies as follows: 5 (top left), 10 (top right), 25 (bottom left), and 50 (bottom right).}
	\label{fig:InitPop}
\end{figure}

The choice of the initial population size $\boldsymbol{N}_0$ is critical for the accuracy of the Gillespie algorithm, with larger values improving the statistical reliability but incurring an increase in memory and computational cost. Moreover, the total population size of the dual process $\boldsymbol{N}_t$ satisfies $\lvert \boldsymbol{N}_t \rvert \sim t^{-1}$, and consequently for sufficiently small $t$ initialising the algorithm with a suitably large $\NN_t$ becomes impractical. To address this issue, we introduce a threshold time $\varepsilon$ and fixed initial population size $N_0$, such that for $t \geq \varepsilon$, the Gillespie algorithm is initialised with total population size $N_0$. We assess the sensitivity of the output to this choice by considering $N_0 \in \{10^4, 10^5, 10^6\}$ and time increments $t \in \{0.01, 0.05, 0.1\}$, and comparing the resulting empirical probability mass functions of $\left|\NN_t\right|$ across all combinations. As shown in Figure~\ref{fig:InitPop}, the empirical distributions coincide for $t \geq 0.05$ across all values of $N_0$, while for $t=0.01$ the choice $N_0=10^4$ exhibits mild discrepancies. Balancing accuracy and computational efficiency, and noting that the algorithm is implemented on a GPU, we set $\varepsilon=0.05$ and $N_0=10^4$. For time increments smaller than $\varepsilon$, a Gaussian approximation is employed, yielding accurate results with manageable runtimes. 
		
\end{document}